# Dwork's conjecture
# on unit root zeta functions

By Daqing Wan*

## 1. Introduction

In this article, we introduce a systematic new method to investigate the conjectural $p$-adic meromorphic continuation of Professor Bernard Dwork's unit root zeta function attached to an ordinary family of algebraic varieties defined over a finite field of characteristic $p$.

After his pioneer $p$-adic investigation of the Weil conjectures on the zeta function of an algebraic variety over a finite field, Dwork went on to study the $p$-adic analytic variation of a family of such zeta functions when the variety moves through an algebraic family. In the course of doing so, he was led to a new zeta function called the unit root zeta function, which goes beyond the reach of the existing theory. He conjectured [8] that such a unit root zeta function is $p$-adic meromorphic everywhere. These unit root zeta functions contain important arithmetic information about a family of algebraic varieties. They are truly $p$-adic in nature and are transcendental functions, sometimes seeming quite mysterious. In fact, no single "nontrivial" example has been proved to be true about this conjecture, other than the "trivial" overconvergent (or $\infty$ log-convergent) case for which Dwork's classical $p$-adic theory already applies; see [6]–[11], [26] for various attempts. In this article, we introduce a systematic new method to study such unit root zeta functions. Our method can be used to prove the conjecture in the case when the involved unit root F-crystal has rank one. In particular, this settles the first "nontrivial" case, the rank one unit root F-crystal coming from the family of higher dimensional Kloosterman sums. Our method further allows us to understand reasonably well about analytic variation of an arithmetic family of such rank one unit root zeta functions, motivated by the Gouvêa-Mazur conjecture about dimension variation of classical and $p$-adic modular forms. We shall introduce another systematic method in a future paper [29] which combined with the method in the present paper will be able to prove Dwork's conjecture in the higher rank case.

---

*This work was partially supported by NSF. The author wishes to thank P. Deligne and N. Katz for valuable discussions at one stage of this work during 1993 and 1994. The author would also like to thank S. Sperber for his careful reading of the manuscript and for his many detailed comments which led to a significant improvement of the exposition.

To explain our results in this introduction, we shall restrict our attention to the family of $L$-functions of higher dimensional Kloosterman sums parametrized by the one-dimensional torus $\mathbf{G}_m$. This is the most concrete example, an essential family for Dwork's conjecture. Let $\mathbf{F}_q$ be the finite field of $q$ elements of characteristic $p$. Let $\Psi$ be a fixed, nontrivial complex-valued, additive character of the finite field $\mathbf{F}_p$. Let $n$ be a positive integer. For each nonzero element $\bar{y} \in \mathbf{F}_q^*$ (the letter $y$ is reserved to denote the Teichmüller lifting of $\bar{y}$) and each positive integer $k$, we define the $n$-dimensional Kloosterman sum over $\mathbf{F}_{q^k}$ to be

$$K_k(\bar{y}, n) = \sum_{x_i \in \mathbf{F}_{q^k}^*} \Psi(\mathrm{Tr}_{\mathbf{F}_{q^k}/\mathbf{F}_p}(x_1 + \cdots + x_n + \frac{\bar{y}}{x_1 \cdots x_n})).$$

The $L$-function attached to the sequence $K_k(\bar{y}, n)$ $(k = 1, 2, \cdots)$ of character sums is defined to be the exponential generating function

$$L(\bar{y}, n, T) = \exp(\sum_{k=1}^{\infty} \frac{K_k(\bar{y}, n)}{k} T^k).$$

It has an Euler product expansion which shows that $L(\bar{y}, n, T)$ is a power series whose coefficients are algebraic integers in the $p^{\mathrm{th}}$ cyclotomic field. The product of all the $p - 1$ conjugates (when $\Psi$ varies over the $p - 1$ nontrivial additive characters of $\mathbf{F}_p$) of $L(\bar{y}, n, T)$ gives the nontrivial part of the zeta function of the affine variety over $\mathbf{F}_q$ defined by

$$z^p - z = x_1 + \cdots + x_n + \frac{\bar{y}}{x_1 \cdots x_n}.$$

Thus, it suffices to understand the above $L$-function.

It is well known that $L(\bar{y}, n, T)^{(-1)^{n-1}}$ is a polynomial of degree $n + 1$. Thus, there are $n + 1$ algebraic integers $\alpha_0(\bar{y}), \cdots, \alpha_n(\bar{y})$ such that

$$L(\bar{y}, n, T)^{(-1)^{n-1}} = (1 - \alpha_0(\bar{y})T) \cdots (1 - \alpha_n(\bar{y})T).$$

Equivalently, for each integer $k \geq 1$, we have the formula

$$(-1)^n K_k(\bar{y}, n) = \alpha_0(\bar{y})^k + \cdots + \alpha_n(\bar{y})^k.$$

By Deligne's theorem [5], these algebraic integers $\alpha_i(\bar{y})$ have complex absolute value $\sqrt{q^n}$. For each prime $\ell \neq p$, these $\alpha_i(\bar{y})$ are $\ell$-adic unit. To describe their $p$-adic absolute values, we fix an embedding of the algebraic numbers $\bar{\mathbf{Q}}$ into $\bar{\mathbf{Q}}_p$ and order the $\alpha_i(\bar{y})$ such that

$$\mathrm{ord}_q \alpha_0(\bar{y}) \leq \cdots \leq \mathrm{ord}_q \alpha_n(\bar{y}).$$

Then Sperber's theorem [21] says that

$$\mathrm{ord}_q\alpha_0(\bar{y}) = 0, \ \mathrm{ord}_q\alpha_1(\bar{y}) = 1, \cdots, \ \mathrm{ord}_q\alpha_n(\bar{y}) = n.$$

This means that there is exactly one characteristic root of slope $i$ for each integer $0 \leq i \leq n$. In particular, there is exactly one $p$-adic unit root $\alpha_0(\bar{y})$ (slope zero).

In a more conceptual language, the $L$-function $L(\bar{y}, n, T)^{(-1)^{n-1}}$ is given by the characteristic polynomial of the Frobenius map at the closed point $\bar{y}$ of a rank $(n+1)$ ordinary overconvergent F-crystal $M_n$ over the torus $\mathbf{G}_m/\mathbf{F}_p$. This F-crystal $M_n$ is not a unit root F-crystal. It is ordinary (Newton polygon coincides with Hodge polygon fibre by fibre). Its unit root part $U_n$ is a sub-F-crystal of rank one on the torus $\mathbf{G}_m/\mathbf{F}_p$, but no longer overconvergent in general. Denote the $1 \times 1$ matrix of the Frobenius map of $U_n$ by $\alpha(y)$. This is a convergent (but not overconvergent in general) power series in $y$ with $p$-adic integral coefficients. The constant term of $\alpha(y)$ is a $p$-adic unit and all other coefficients of $\alpha(y)$ are divisible by $\pi$, where $\pi^{p-1} = -p$. Let $R$ be the ring of integers in $\mathbf{Q}_p(\pi)$. Then, we can write

$$\alpha(y) = \sum_{i=0}^{\infty} b_i y^i, \ b_0 \in R^*, \ b_i \in \pi R \ (i > 0), \ \lim_{i \to \infty} b_i = 0.$$

This power series $\alpha(y)$ can be expressed explicitly in terms of a certain uniquely determined solution of the hypergeometric differential equation attached to the family of Kloosterman sums. If a nonzero element $\bar{y} \in \mathbf{F}_q^*$ with $q = p^a$, then the unit root $\alpha_0(\bar{y})$ of the above $L$-function of Kloosterman sums over $\mathbf{F}_q$ is given by the $p$-adic analytic formula

$$\alpha_0(\bar{y}) = \alpha(y)\alpha(y^p)\cdots\alpha(y^{p^{a-1}}),$$

where $y$ is the Teichmüller lifting of $\bar{y}$. This formula is a consequence of the Hodge-Newton decomposition discovered by Dwork [7], [8] in the early seventies.

To further understand how the unit root $\alpha_0(\bar{y})$ (as a $p$-adic integer) varies when $\bar{y}$ varies, one naturally introduces and considers the $L$-function $L(U_n, T)$ attached to the unit root F-crystal $U_n$. More generally, for an integer $k$, we can consider the $k^{\mathrm{th}}$ tensor power $U_n^{\otimes k}$. Its $L$-function over the prime field $\mathbf{F}_p$ is defined in the usual way by the infinite Euler product

$$L(U_n^{\otimes k}, T) = \prod_{\bar{y} \in \mathbf{G}_m} \frac{1}{(1 - T^{\deg(\bar{y})}\alpha_0^k(\bar{y}))}$$

$$= \prod_{\bar{y} \in \mathbf{G}_m} \frac{1}{(1 - T^{\deg(\bar{y})}\alpha^k(y)\alpha^k(y^p)\cdots\alpha^k(y^{p^{\deg(\bar{y})-1}}))},$$

where $\bar{y}$ runs over the closed points of $\mathbf{G}_m/\mathbf{F}_p$. This $L$-function $L(U_n^{\otimes k}, T)$ can also be defined using the exponential generating function associated to the sequence of $p$-adic character sums arising from the $p$-adic character $\alpha(y)$. The unit root F-crystal $U_n$ gives a one-dimensional continuous $p$-adic representation $\rho_n$ of the arithmetic fundamental group $\pi_1^{\mathrm{arith}}(\mathbf{G}_m/\mathbf{F}_p)$:

$$\rho_n : \pi_1^{\mathrm{arith}}(\mathbf{G}_m/\mathbf{F}_p) \longrightarrow \mathrm{GL}_1(R),$$

where the torus $\mathbf{G}_m/\mathbf{F}_p$ is the parameter space for the parameter $\bar{y}$. The $L$-function $L(U_n^{\otimes k}, T)$ is the same as the $L$-function $L(\rho_n^{\otimes k}, T)$ attached to the continuous $p$-adic representation $\rho_n^{\otimes k}$. According to a plausible conjecture of Katz [14] (see [11] for a proof in the constant sheaf case), this mysterious unit root zeta function $L(U_n, T)$ ($k = 1$) is also the $L$-function attached to the relative $p$-adic étale cohomology with compact support of the family of Kloosterman sheaves parametrized by the torus $\mathbf{G}_m$. These are transcendental, nonrational functions. The meromorphic continuation of $L(U_n^k, T)$ to the closed unit disc already implies the existence of a weak $p$-adic equi-distribution theorem for the $p$-adic "angle" of the zeros of the $L$-function $L(\bar{y}, n, T)$ of Kloosterman sums [17]. Dwork's unit root conjecture [8] is the following:

CONJECTURE (DWORK).    *For every integer $k$, the unit root zeta function $L(U_n^{\otimes k}, T)$ is $p$-adic meromorphic.*

For a so-called overconvergent F-crystal, the $L$-function is always meromorphic by Dwork's trace formula. The difficulty about this conjecture is that the unit root F-crystal $U_n$ (obtained by solving the fixed point of a contraction map in a $p$-adic Banach space) is no longer overconvergent in general. In the case that $n = 1$, Dwork [9] showed that there is an excellent lifting (Deligne-Tate lifting) of the Frobenius map $x \to x^p$ such that $U_1$ is overconvergent with respect to the excellent lifting. Thus, the classical overconvergent theory applies and the zeta function $L(U_1^k, T)$ is $p$-adic meromorphic in the special case $n = 1$.

Unfortunately, the so-called excellent lifting rarely exists. In fact, for each $n > 1$, Sperber [21] showed that there does not exist an excellent lifting of the Frobenius map such that $U_n$ is overconvergent. Thus, the situation cannot be reduced to the "trivial" overconvergent case. Although it was conjectured more generally by Katz [14] that the $L$-function of any F-crystal is always $p$-adic meromorphic, this was disproved in [26]. Thus, to handle Dwork's conjecture for $n > 1$, one has to employ a new approach. The purpose of this article is to introduce a new method which can be used to handle the case when the unit root F-crystal is of rank one. See Section 3 for a full description of our results. In particular, we obtain the following:

THEOREM 1.1. *For every integer $k$, the $k^{\text{th}}$ unit root zeta function $L(U_n^{\otimes k}, T)$ is $p$-adic meromorphic.*

The general tool for $p$-adic meromorphic continuation of $L$-functions is to use Dwork's trace formula. It expresses the unit root zeta function as an alternating product of the Fredholm determinants of several continuous operators in a $p$-adic Banach space. If these operators were completely continuous [20] (which is the case if the unit root F-crystal is overconvergent or more generally $\infty$ log-convergent), then the Fredholm determinants would be $p$-adic entire and hence the $L$-function would be meromorphic. Unfortunately, those operators do not seem to be completely continuous in general. All previous approaches to Dwork's conjecture try to prove that the Fredholm determinants are entire. This is usually done by showing that the radius of convergence is infinite. But it may well be the case that those Fredholm determinants could be meromorphic but *not* entire. And so, the radius of convergence of those Fredholm determinants could indeed be finite because of possible poles. For an arbitrary unit root F-crystal, the Fredholm determinant is not even meromorphic, not mentioning entire, as the counterexample in [26] shows.

Our new approach takes a long detour and avoids proving the likely false statement that the Fredholm determinants are entire. The basic idea is to express each $L(U_n^{\otimes k}, T)$ as the limit of a sequence of functions which are alternating products of $L$-functions of various higher symmetric powers and exterior powers of the original overconvergent F-crystal $M_n$ we started with. These later $L$-functions are meromorphic since their F-crystals are overconvergent. However, they become out of control if we try to take the limit. In the case that the unit root F-crystal $U_n$ is of rank one (which is always the case for our Kloosterman family), we are able to prove that our sequence of meromorphic $L$-functions are uniformly meromorphic. This uniformity and a continuity result allow us to take the limit and deduce that the limiting unit root zeta function $L(U_n^{\otimes k}, T)$ is still a meromorphic function. The actual proof decomposes $L(U_n^{\otimes k}, T)$ as a quotient of two functions, each of which is the limit of a sequence of entire functions. This decomposition makes it easier to take the limit. Note that our work does not imply that the above mentioned Fredholm determinants are entire. It only implies that they are meromorphic which is all one needs to prove Dwork's conjecture.

Once we know that $L(U_n^{\otimes k}, T)$ is meromorphic for each $k$, we could ask how the unit root zeta function $L(U_n^{\otimes k}, T)$ varies with the integer parameter $k$. This can be viewed as an extension of Gouvêa-Mazur's conjecture [12], [13] about dimension variation of modular forms. The Gouvêa-Mazur conjecture corresponds exactly to the case of unit root zeta functions arising from the ordinary family of elliptic curves over a finite field, in which case the unit root F-crystal is of rank 1 and already overconvergent by Deligne-Tate lifting. The unit root zeta function in this elliptic case, more precisely the special value $L(U_n^{\otimes k}, 1)$, is also directly related to geometric Iwasawa theory and $p$-adic $L$-functions, see [3], [4], [18]. Our results here are strong enough to obtain useful information about analytic variation of the family of unit root zeta functions $L(U_n^{\otimes k}, T)$ as the integer parameter $k$ varies.

The definition of $L(U_n^{\otimes k}, T)$ and Fermat's little theorem show that the $L$-function $L(U_n^{\otimes k}, T)$ is $p$-adically continuous in $k$. More precisely, if $k_1$ and $k_2$ are two integers such that

$$(1.1) \qquad k_1 \equiv k_2 \ (\mathrm{mod} \ (p-1)p^m),$$

then

$$(1.2) \qquad L(U_n^{\otimes k_1}, T) \equiv L(U_n^{\otimes k_2}, T) \ (\mathrm{mod} \ \pi p^m).$$

To get deeper information about the zeros of $L(U_n^{\otimes k}, T)$, we would like to decompose $L(U_n^{\otimes k}, T)$ in terms of its slopes.

For a given integer $k$, the $p$-adic meromorphic function $L(U_n^{\otimes k}, T)$ can be factored completely over $\bar{\mathbf{Q}}_p$:

$$L(U_n^{\otimes k}, T) = \prod_{i \geq 0}(1 - z_i(k)T)^{\pm 1},$$

where each $z_i(k)$ is a $p$-adic integer in $\bar{\mathbf{Q}}_p$. For a rational number $s \in \mathbf{Q}$, we define the slope $s$ part of $L(U_n^{\otimes k}, T)$ to be

$$L_s(U_n^{\otimes k}, T) = \prod_{\mathrm{ord}_p z_i(k) = s}(1 - z_i(k)T)^{\pm 1}.$$

This is a rational function. Let $d_s(U_n, k)$ denote the degree of the slope $s$ rational function $L_s(U_n^{\otimes k}, T)$, where the degree of a rational function means the degree of its numerator minus the degree of its denominator. This function $d_s(U_n, k)$ of two variables is called the degree function of the meromorphic $L$-function $L(U_n^{\otimes k}, T)$. Note that $d_s(U_n, k)$ can take negative values. A fundamental problem is to understand the degree function $d_s(U_n, k)$; see [28] for a discussion about possible upper bounds for $d_s(U_n, k)$ in the elliptic case. Here, following the spirit of the Gouvêa-Mazur conjecture, we only investigate how the function $d_s(U_n, k)$ varies with $k$ when $s$ is fixed.

For any rational number $s$, we define $m(U_n, s)$ to be the smallest nonnegative integer $m \in \mathbf{Z}_{\geq 0} \cup \{+\infty\}$ such that whenever $k_1$ and $k_2$ are two integers satisfying (1.1), we have the equality

$$d_{s^*}(U_n, k_1) = d_{s^*}(U_n, k_2)$$

for all $0 \leq s^* \leq s$. By our definition, the following trivial inequality

$$0 \leq m(U_n, s) \leq +\infty$$

holds for each $s$. Furthermore, $m(U_n, s) = 0$ for $s < 0$. This is because $d_s(U_n, k) = 0$ for $s < 0$. The quantity $m(U_n, s)$ gives information about how often the degree function $d_s(U_n, k)$ changes as $k$ varies $p$-adically, where the slope $s$ is fixed. Our result implies:

THEOREM 1.2. *The number $m(U_n, s)$ is finite for every $s$ and every $n$.*

This gives a control theorem like the one Coleman [3] obtained for the family of elliptic curves in his work on the Gouvêa-Mazur conjecture about dimension variation of modular forms. In particular, it implies that for each fixed $s$, $d_s(U_n, k)$ is a $p$-adically locally constant function of $k$ if we restrict $k$ to a residue class modulo $p - 1$. Thus, the degree function $d_s(U_n, k)$ is a bounded function of $k$ for each fixed $s$. Our work makes it possible and meaningful to extend the Gouvêa-Mazur conjecture from the family of elliptic curves to a family of algebraic varieties with one unit root, fibre by fibre. Since our method is effective, it can be used to give an explicit upper bound for $m(U_n, s)$, generalizing the quadratic bound in [27] about Gouvêa-Mazur's conjecture. We shall, however, not make our calculations explicit in this paper because we work under the weakest possible condition of a $c$ log-convergent setting, where the control theorem is already best possible in some sense. For those unit root F-crystals which are embedded in overconvergent F-crystals (such as Kloosterman sums), there is a somewhat simpler version of our method which would give a reasonably good and explicit estimate depending on several factors. In fact, in a future article on overconvergent setup, we will prove a general result which implies immediately that the quantity $m(U_n, s)$ is bounded by a polynomial in $s$ of degree $n + 2$. In the special case $n = 1$, this is a cubic bound which can be improved to be quadratic because of the existence of an excellent lifting.

An important further question is to ask about the possible slopes or the $p$-adic absolute values of the zeros and poles of the unit root zeta function $L(U_n^{\otimes k}, T)$, namely the $p$-adic Riemann hypothesis of the meromorphic function $L(U_n^{\otimes k}, T)$. Let $S(k, n)$ be the set of slopes of the zeros and poles of the $L$-function $L(U_n^{\otimes k}, T)$, where the slope of a $p$-adic number $\alpha$ means $\text{ord}_p \alpha$. This is in general an infinite set of rational numbers. The $p$-adic Riemann hypothesis in this case is the following conjecture. We stay with the simpler $n = 1$ case since so little is known.

CONJECTURE 1.3. *Let $n = 1$. Then, the set $S(k, 1)$ has a uniformly bounded denominator for all integers $k$.*

This conjecture is not known to be true for a single $k$ other than the exceptional value $k = 0$ in which case

$$L(U_n^{\otimes 0}, T) = L(1, T) = (1 - T)/(1 - pT)$$

is rational. A slightly stronger version of Conjecture 1.3 is the statement that the two variable degree function $d_s(U_1, k)$ is uniformly bounded for $n = 1$. One may wonder why we put the restriction $n = 1$ in the above conjecture. The main reason is that for $n = 1$, we can at least prove the weak evidence that a suitable average version of Conjecture 1.3 is true by using the existence of an excellent lifting. For $n > 1$, we are currently unable to prove any probabilistic version of Conjecture 1.3. Thus, it seems safer to stay with the simplest case $n = 1$ at the moment.

Our results extend to a higher slope case (again assuming rank 1) as also conjectured by Dwork. The extension to higher slopes uses more tensor products. The first step is to separate the higher slope unit root zeta function from other smaller slope unit root zeta functions in a natural way. After the separation is done, it turns out that the higher slope unit root zeta function can also be expressed, in a more complicated way, as the limit of a sequence of meromorphic functions. The necessary uniformity and continuity result as the slope zero case can be established when we assume the concerned unit root F-crystal is of rank one. This condition is satisfied for the Kloosterman family by Sperber's theorem.

If the unit root F-crystal has rank greater than one, the uniform part of our approach does not work. There are easy examples (constant unit root F-crystals of rank greater than one, for instance) which show that the sequence of the $L$-functions of the higher symmetric powers are not uniformly meromorphic and thus one cannot pass the meromorphy property to its limit. Although one can use the results of this paper to prove Dwork's conjecture in some higher rank cases, we shall not pursue it here as it does not seem to get very far. We shall introduce another method in a future article [29] which, combined with the method in the present paper will be able to prove Dwork's conjecture in full generality.

For simplicity, we shall work in the toric setting in this paper so that Dwork's trace formula can be applied and so that we can go as far as possible (proving optimal results) without too much technical burden. For a general smooth base scheme, one needs to use Reich-Monsky's generalization of Dwork's trace formula. It would be more complicated and involve more nasty analytic arguments. We shall also forget the horizontal connection (differential equation) of an F-crystal and work in a larger category called $\sigma$-modules of possibly infinite rank with certain nuclear structure. Thus, our results actually go beyond, in several ways, the conjectural context of F-crystals of finite rank. In particular, our overconvergent nuclear $\sigma$-module (even of finite rank) is generally not overconvergent in Berthelot's sense [1] if one wants to consider it as an F-crystal by adding back the implied horizontal connection of the Frobenius map. The difference is that we only assume that the Frobenius map is overconvergent. Berthelot's rigid cohomology assumes that both the Frobenius map and the horizontal connection are overconvergent. For example, the rank one unit root F-crystal coming from the ordinary family of elliptic curves is overconvergent in our sense but not overconvergent in Berthelot's sense.

As indicated in [25], there is a characteristic $p$ version of Dwork's conjecture which would be useful in studying analytic variation of $L$-functions of $\varphi$-sheaves and Drinfeld modules [22]. Our method here shows that the characteristic $p$ Dwork conjecture holds in the rank one case as well, where the base scheme is allowed to be an arbitrary separated scheme of finite type over $\mathbf{F}_q$. In the characteristic $p$ situation, the general base scheme case can be reduced to the toric case by a simple trick because the Frobenius lifting problem disappears in characteristic $p$. Although the characteristic $p$ version of Dwork's conjecture is technically simpler, it requires the same amount of original ideas. In particular, we do not know how to use the method of the present paper to prove the characteristic $p$ Dwork conjecture if the unit root part has rank greater than one.

From a structural point of view, our result simply says that Dwork's unit root zeta function attached to the rank one unit root part of an ordinary overconvergent F-crystal of finite rank is a finite alternating product of $L$-functions of certain overconvergent "F-crystals of infinite rank". This gives an explicit formula for Dwork's unit root zeta function in terms of $L$-functions of overconvergent "F-crystals of infinite rank". We shall show [29] that a more complicated formula holds for a higher rank unit root zeta function, but we will need infinitely many overconvergent "F-crystals of infinite rank". To handle these new objects (overconvergent "F-crystals of infinite rank") called overconvergent nuclear $\sigma$-modules in this paper, one often needs to be careful about various convergent conditions. The nuclear condition (some sort of infinite Hodge structure) can be viewed as a weak finiteness condition. Our nuclear operators are not completely continuous in the usual sense (in fact our Banach space does not even have an orthonormal basis). However, they can be viewed as the dual of completely continuous operators. Working in the dual noncompletely continuous but nuclear category has some significant advantages in our proofs. Our nuclear condition insures that everything works as expected. There is another reason which will be made more explicit in a future paper. Namely, even if one starts with an infinite rank $\sigma$-module for the more familiar Banach space with an orthonormal basis, in extending the more abstract Monsky trace formula to the infinite rank situation, one would encounter a suitable dual which would land in the category of Banach spaces without an orthonormal basis. Thus, it makes sense to deal with them at the beginning.

The setting of this paper starts immediately with the infinite rank case. However, those who are interested only in Dwork's original conjecture can assume that all our objects are in the familiar finite rank case. This is enough to prove the rank one case of Dwork's conjecture in his original finite rank setting. To study further analytic variation of unit root zeta functions, one can still stay in the finite rank case but it is a little awkward to do it because one has to take the limit of a discontinuous family. We get around the difficulty by working with "essentially continuous families". Thus, for the analytic variation part, we actually have two proofs. One uses the notion of "essential continuity" and stays in the finite rank case. The other uses continuous families but works in the infinite rank setting. The infinite rank setting also gives us the flexibility to get around the difficulty caused by the (unknown) conjectural generic finite dimensionality of the relative rigid cohomology of a family of varieties $f : X \to Y$ over a finite field. This generic finite dimensionality is needed to reduce the geometric version of Dwork's conjecture to the more general version about $L$-functions of F-crystals if one insists on staying in the finite rank case. The conjectural finite dimensionality has recently been proved by Berthelot [2] when the parameter space $Y$ is a point. The generic finite dimensionality of the relative rigid cohomology in general is still open. Ultimately, we believe that one should work with the infinite rank setting because that is essentially where the unit root zeta function sits and that is where the unit root zeta function can be best understood. For these reasons, we have chosen to work directly with the infinite (countable) rank setting throughout the paper.

In addition to the example of the higher dimensional Kloosterman family described in this introduction, we note that there are many other families which have exactly one $p$-adic unit root fibre by fibre. For instance, if $f(x, y)$ is an arbitrary Laurent polynomial over $\mathbf{F}_q$ in two sets of variables $x = (x_1, \cdots, x_n)$ and $y = (y_1, \cdots, y_m)$, then we can view $f(x, y)$ as a family of Laurent polynomials in $x$ parametrized by the parameter $y$ varying in the $m$-dimensional torus. For each closed point $y$ in the $m$-torus, we can form a sequence of exponential sums of the Laurent polynomial $f(x, y)$ on the $n$-torus and hence get an $L$-function $L(f(x, y), T)$ for each closed point $y$ of the $m$-torus. The family of $L$-functions $L(f(x, y), T)^{(-1)^{n-1}}$ parametrized by $y$ of the $m$-torus has exactly one $p$-adic unit root fibre by fibre. Our result applies to this much more general situation. In the $n$-dimensional Kloosterman family discussed above, the Laurent polynomial $f(x, y)$ is simply given by $x_1 + \cdots + x_n + y/x_1 \cdots x_n$ with $m = 1$.

The content is organized as follows. Section 2 introduces the basic notion of a nuclear $\sigma$-module and defines its $L$-function. This can be viewed as an extension of F-crystals without connection from finite rank to infinite rank.

Section 3 reviews the Hodge-Newton decomposition of an ordinary F-crystal in the context of nuclear $\sigma$-modules. This allows us to reformulate Dwork's conjecture in a more general framework where adequate tools can be developed. Section 4 gives a basic decomposition formula for Dwork's $k^{\text{th}}$ power $L$-function of a nuclear $\sigma$-module in terms of more natural $L$-functions of various higher symmetric powers and exterior products. This formula is our starting point. It already gives a nontrivial result about meromorphic continuation of Dwork's unit root zeta function. Section 5 studies continuous and uniform families of $p$-adic meromorphic functions. It is in this section that various basic uniform results are proved. Section 6 begins the proof of the rank one slope zero part of Dwork's conjecture by taking limit and showing various continuity results. Section 7 extends our slope zero result to the higher slope case (still rank one). Section 8 describes an alternating approach to the limiting procedure. It includes an explicit formula for Dwork's unit root zeta function in the rank one case as well as a slightly more general rank one result which is crucial in our forthcoming work on the higher rank case of Dwork's conjecture.

## 2. Nuclear $\sigma$-modules and $L$-functions

In the introduction, the base field for the unit root zeta function is the prime field $\mathbf{F}_p$. From now on, we will work with the slightly more general base field $\mathbf{F}_q$. Let $K$ be a fixed finite extension of the field $\mathbf{Q}_p$ of $p$-adic rational numbers with residue field $\mathbf{F}_q$. Let $\pi$ be a fixed uniformizer of $K$, and $R$ be the ring of integers in $K$. This ring $R$ consists of those elements $a \in K$ such that $\operatorname{ord}_\pi a \geq 0$. We normalize the $p$-adic ($\pi$-adic) absolute value on $K$ by defining $|\pi|_\pi = 1/p$. Let $n$ be a fixed positive integer. We shall consider various nuclear $\sigma$-modules over the $n$-dimensional torus $\mathbf{G}_m^n$, which are certain Banach modules with a nuclear action of a semi-linear operator. Throughout the paper, our base space will be the $n$-dimensional torus $\mathbf{G}_m^n$.

First, we define various coefficient rings. Let

$$A_0 = \{ \sum_{u \in \mathbf{Z}^n} a_u X^u | a_u \in R, \ \lim_{|u| \to \infty} a_u = 0 \}$$

be the ring of *convergent* Laurent series over $R$, where for a lattice point $u = (u_1, \cdots, u_n) \in \mathbf{Z}^n$ we define $X^u = X_1^{u_1} \cdots X_n^{u_n}$ and $|u| = |u_1| + \cdots + |u_n|$. Let $\sigma$ be the $R$-linear Frobenius map on $A_0$ defined by:

$$\sigma(X^u) = X^{qu}.$$

Thus, one can think of the pair $(A_0, \sigma)$ as a natural $p$-adic lifting of the characteristic $p$ coordinate ring of the $n$-torus $(\mathbf{G}_m^n, \bar{\sigma})$ over $\mathbf{F}_q$, where $\bar{\sigma}$ is the $q^{\text{th}}$ power Frobenius map. The ring $A_0$ is $p$-adically complete. However, it is too large for $L$-functions and cohomological purpose. We shall frequently need to work in various subrings of $A_0$ which are not complete but weakly complete in some sense in order to get interesting results about $L$-functions.

The *overconvergent* subring $A^\dagger$ of $A_0$ consists of those elements satisfying

$$\lim_{|u| \to \infty} \inf \frac{\operatorname{ord}_\pi a_u}{|u|} > 0.$$

For a real number $0 \leq c \leq \infty$, let $A_c$ be the $c \log$-*convergent* subring of $A_0$ consisting of those elements satisfying

$$\lim_{|u| \to \infty} \inf \frac{\operatorname{ord}_\pi a_u}{\log_q |u|} \geq c.$$

It is clear that for $c_1 > c_2 > 0$, we have the inclusion relation

$$R[X] \subset A^\dagger \subset A_\infty \subset A_{c_1} \subset A_{c_2} \subset A_0.$$

All of these rings are known to be Noetherian [24]. They are not $p$-adically complete except for the largest one $A_0$. The algebra $A_0$ is a complete normed algebra under the Gauss norm:

$$\| \sum_u a_u X^u \| = \max_u |a_u|_\pi.$$

Thus, the ring $A_0$ becomes a Banach algebra over $R$ and the monomials $\{X^u\}$ form an orthonormal basis of $A_0$ over $R$. A Banach module $M$ over $A_0$ is an ultra normed complete $A_0$-module, whose norm is also denoted by $\| \ \|$. It satisfies the relations

$$\|m\| = 0 \text{ if and only if } m = 0,$$

$$\|am\| = \|a\|\|m\|, \text{ for } a \in A_0, m \in M,$$

$$\|m_1 + m_2\| \leq \max(\|m_1\|, \|m_2\|), \text{ for } m_1, m_2 \in M.$$

A (formal) *basis* of $M$ over $A_0$ is a subset $\{e_i : i \in I\}$ of $M$ for some index set $I$, such that every element in $M$ can be written uniquely in the form $\sum a_i e_i$ with $a_i \in A_0$. Thus, a Banach module $M$ over $A_0$ with a basis is the space of sequences $\{a_i, i \in I\}$ with entries in $A_0$. An *orthonormal basis* of $M$ over $A_0$ is a subset $\{e_i : i \in I\}$ of $M$ for some index set $I$, such that every element in $M$ can be written uniquely in the form $\sum a_i e_i$ with $a_i \in A_0$ and the additional condition $\lim_i \|a_i\| = 0$. Thus, a Banach module $M$ over $A_0$ with an orthonormal basis is the space of convergent sequences with entries in $A_0$. Note that a Banach module with a basis does not necessarily have an orthonormal basis. For example, the formal power series ring $R[[x]]$ is a Banach space with the basis $\{x^u\}$ over $R$ but it does not have an orthonormal basis over $R$.

For a finite rank free module $M$ over $A_0$, basis and orthonormal basis are the same concept. Any free $A_0$-module $M$ of finite rank can be viewed as a Banach module of finite rank by suitably extending the norm on $A_0$ to $M$, depending on a chosen basis $\vec{e}$. We remind the reader that our definition of a (formal) basis in the infinite rank case may be different from the traditional usage of the notion "basis", where it is assumed that the sum $\sum a_i e_i$ is a finite combination. We do allow formal infinite combination. The word "basis" in this paper will always refer to formal basis. We shall be interested in the Banach module over $A_0$ of the form

$$M = \{\sum_i a_i e_i | a_i \in A_0\},$$

where $\vec{e} = \{e_1, \cdots, e_i, \cdots\}$ is a basis indexed by a set $I$ of at most countable cardinality (rank). We assume that the norm on $M$ is defined by

$$\|\sum_i a_i e_i\| = \max_i \|a_i\|.$$

If the rank is infinite, this Banach module $M$ has a basis over $A_0$ but does not have an orthonormal basis. If $(b_{i,j})$ is a matrix with entries $b_{i,j}$ in $A_0$ or $R$, we use $\mathrm{ord}_\pi(b_{i,j})$ to denote $\min_{i,j} \mathrm{ord}_\pi b_{i,j}$. Let $c(I)$ be the set of convergent sequences over $A_0$ indexed by $I$. Let $b(I)$ be the set of sequences over $A_0$ indexed by $I$. The index set $I$ will be identified with a subset of the positive integers.

We first study which transition matrix gives a new basis for $M$ over $A_0$.

LEMMA 2.1. *Let $\vec{f} = \vec{e}U$ be a set of elements in $M$ indexed by $I$, where $U$ is a matrix whose rows and columns are both indexed by $I$. Then, $\vec{f}$ is a basis of $M$ over $A_0$ if and only if $U$ is invertible over $A_0$ and the row vectors of both $U$ and $U^{-1}$ are in $c(I)$. In this case, there exists the norm relation*

$$\|\sum_i a_i f_i\| = \max_i \|a_i\|.$$

*Proof.* Write

$$f_j = \sum_i u_{ij} e_i.$$

If $\vec{f}$ is a basis of $M$ over $A_0$, then $U$ is invertible and the infinite sum

$$\sum_j f_j = \sum_i (\sum_j u_{ij}) e_i$$

is convergent. Thus, the row vectors of $U$ are in $c(I)$. The same statement holds for $U^{-1}$ since $\vec{e} = \vec{f}U^{-1}$ is a basis.

Conversely, assume that the row vectors of $U$ and $U^{-1}$ are in $c(I)$. Then, for any column vector $\vec{a} \in b(I)$, the substitution

$$\vec{e}\,\vec{a} = \vec{f}(U^{-1}\vec{a})$$

is well-defined. This means that each element of $M$ can be written as a combination of the elements in $\vec{f}$. If we have two such expressions of the same element:

$$\vec{f}\,\vec{a} = \vec{f}\,\vec{b},$$

the substitution

$$\vec{e}(U\vec{a}) = \vec{e}(U\vec{b})$$

is well-defined. Since $\vec{e}$ is a basis, we must have

$$U\vec{a} = U\vec{b}.$$

Since the row vectors of $U^{-1}$ are in $c(I)$, we can multiply the above equation by $U^{-1}$ on the left side. This implies that $\vec{a} = \vec{b}$. Thus, $\vec{f}$ is a basis.

To prove the norm relation, it suffices to prove for a column vector $\vec{a} \in b(I)$ not divisible by $\pi$, we have $\|\vec{f}\,\vec{a}\| = \|\vec{e}(U\vec{a})\| = 1$. Since $U$ is invertible and $\vec{a}$ is not divisible by $\pi$, it follows that $U\vec{a}$ is not divisible by $\pi$. This implies that $\|\vec{e}(U\vec{a})\| = 1$. The proof is complete.                    □

*Definition 2.2.* *A nuclear $\sigma$-module over $A_0$ is a pair $(M, \phi)$, where $M$ is a Banach module over $A_0$ with a basis $\vec{e} = (e_1, e_2, \cdots)$ indexed by $I$ of at most countable cardinality such that $\|\sum_i a_i e_i\| = \max_i \|a_i\|$, and $\phi$ is an infinite $\sigma$-linear map:*

$$\phi : M \longrightarrow M, \ \phi(\sum_{i \geq 1} a_i e_i) = \sum_{i \geq 1} \sigma(a_i)\phi(e_i), \ a_i \in A_0$$

*satisfying the nuclear condition*

$$(2.1) \qquad\qquad\qquad \lim_{i \to \infty} \|\phi(e_i)\| = 0.$$

Note that the map $\phi$ in Definition 2.2 is automatically continuous since $\phi$ is bounded. In fact, $\|\phi(m)\| \leq 1$ for all $m \in M$. We emphasize that our Banach module $M$ has a basis but does not have an orthonormal basis in the infinite rank case. Thus, specifying the image of $\phi$ at the basis elements $e_i$ is not enough to define an infinite $\sigma$-linear map. This is because the infinite sum $\phi(\sum a_i e_i) = \sum_i \sigma(a_i)\phi(e_i)$ may not be defined when we write out $\phi(e_i)$. Our condition that $\lim_i \|\phi(e_i)\| = 0$ insures that the infinite sum $\sum_i \sigma(a_i)\phi(e_i)$ is indeed well-defined. Although our nuclear $\sigma$-module is defined in terms of a chosen basis, the resulting concept is actually independent of the chosen basis. Now, we have the following:

LEMMA 2.3.    *If $(M, \phi)$ is nuclear with respect to one basis $\vec{e}$, then $(M, \phi)$ is nuclear with respect to any other basis $\vec{f}$.*

*Proof.* We need to show that

$$(2.2) \qquad\qquad\qquad \lim_{i \to \infty} \|\phi(f_i)\| = 0.$$

Write

$$f_j = \sum_i u_{ij} e_i.$$

Since $\phi$ is nuclear with respect to $\vec{e}$, for any $\epsilon > 0$, there is an integer $N > 0$ such that for all integers $i > N$, we have $\|\phi(e_i)\| < \epsilon$. This implies that

$$\|\phi(f_j)\| = \|\sum_i u_{ij}^\sigma \phi(e_i)\|$$

$$\leq \max(\epsilon, \|u_{1j}^\sigma \phi(e_1)\|, \cdots, \|u_{Nj}^\sigma \phi(e_N)\|)$$

$$\leq \max(\epsilon, \|u_{1j}\|, \cdots, \|u_{Nj}\|).$$

Since for each fixed $i$, we have by Lemma 2.1 that

$$\lim_j \|u_{ij}\| = 0.$$

There is then an integer $N_1 > 0$ such that for all integers $1 \leq i \leq N$ and all integers $j > N_1$, we have $\|u_{ij}\| < \epsilon$. This proves that for all $j > \max(N_1, N)$, we have

$$\|\phi(f_j)\| < \epsilon.$$

The lemma is proved. □

If $M$ is of finite rank, the nuclear condition (2.1) is automatically satisfied. Thus, in the finite rank case, we may drop the word "nuclear" and simply say a $\sigma$-module. Intuitively, one could think of our nuclear map $\phi$ as the dual of a family of completely continuous semi-linear operators (parametrized by the variable $X$) with an orthonormal basis. In fact, if $B(X)$ is the matrix of $\phi$ under a basis, then the transpose of $B(X)$ can be viewed as the matrix of a family of completely continuous operators of a Banach space with an orthonormal basis. Thus, the standard $p$-adic spectral theory applies fibre by fibre, although the fibres are generally not completely continuous in the usual sense.

A morphism between two nuclear $\sigma$-modules $(M, \phi)$ and $(N, \psi)$ is an $A_0$-linear map of an $A_0$-module

$$\theta : M \longrightarrow N$$

such that $\theta \circ \phi = \psi \circ \theta$. In this way, the category of nuclear $\sigma$-modules is defined. In particular, it makes sense to talk about isomorphic nuclear $\sigma$-modules. It is easy to check that the usual direct sum of two nuclear $\sigma$-modules

$$(M, \phi) \oplus (N, \psi) = (M \oplus N, \phi \oplus \psi)$$

is again a nuclear $\sigma$-module. To extend the usual tensor operations to nuclear $\sigma$-modules, we need to introduce the concepts of formal tensor product, formal symmetric powers and formal exterior powers. We shall restrict to the definition of formal tensor product since the other concepts are similar.

*Definition* 2.4. Let $(M, \phi)$ be a nuclear $\sigma$-module with a basis $\vec{e}$. Let $(N, \psi)$ be a nuclear $\sigma$-module with a basis $\vec{f}$. Let

$$\phi(e_i) = \sum_k u_{ki} e_i, \ \psi(f_j) = \sum_k v_{kj} f_j.$$

Define the formal tensor product of $(M, \phi)$ and $(N, \psi)$ to be the nuclear $\sigma$-module $(M \otimes N, \phi \otimes \psi)$, where $M \otimes N$ denotes the Banach $A_0$-module with the basis

$$\vec{e} \otimes \vec{f} = (\cdots, e_i \otimes f_j, \cdots)$$

and $\phi \otimes \psi$ is determined by the relation

$$
\begin{aligned}
\phi \otimes \psi(e_i \otimes f_j) &= \phi(e_i) \otimes \psi(f_j) \\
&= (\sum_{k_1} u_{k_1 i} e_{k_1}) \otimes (\sum_{k_2 j} v_{k_2 j} f_{k_2}) \\
&= \sum_{k_1, k_2} u_{k_1 i} v_{k_2 j} (e_{k_1} \otimes f_{k_2}).
\end{aligned}
$$

It is easy to check that

$$\lim_{i+j \to \infty} \|\phi \otimes \psi(e_i \otimes f_j)\| = 0.$$

Thus, the *formal tensor product* $(M \otimes N, \phi \otimes \psi)$ is indeed a *nuclear $\sigma$-module*.

Note that we are using the usual tensor notation for the formal tensor product. In fact, we shall often call a formal tensor product just a tensor product. This should not cause confusion since all tensor products will be formal tensor products in this paper. In the finite rank case, of course, the formal tensor product agrees with the usual tensor product. Similarly, one can define formal symmetric powers $(\mathrm{Sym}^k M, \mathrm{Sym}^k \phi)$ and formal exterior powers $(\wedge^k M, \wedge^k \phi)$, using the same classical notation. All symmetric powers and all exterior powers in this paper are in the formal sense. Our definition of the formal tensor product $(M \otimes N, \phi \otimes \psi)$ is based on a basis $\vec{e}$ of $M$ and a basis $\vec{f}$ of $N$. The resulting concept is actually independent of the choice of the bases.

LEMMA 2.5. *The formal tensor product of $(M, \phi)$ and $(N, \psi)$ is independent of the choice of the basis $\vec{e}$ of $M$ and the basis $\vec{f}$ of $N$.*

*Proof.* Let $\vec{e^*}$ (resp. $\vec{f^*}$) be another basis of $M$ (resp. $N$). Write

$$\vec{e^*} = \vec{e} U, \ \lim_{j \to \infty} \|u_{ij}\| = 0,$$

$$\vec{f^*} = \vec{e} V, \ \lim_{j \to \infty} \|v_{ij}\| = 0.$$

We use Lemma 2.1 to show that $\vec{e^*} \otimes \vec{f^*}$ is a basis of the formal tensor product $M \otimes N$. The transition matrix $U \otimes V$ between $\vec{e} \otimes \vec{f}$ and $\vec{e^*} \otimes \vec{f^*}$ is clearly invertible. We need to check the convergent property of the row vectors in the transition matrix. Since

$$e_i^* \otimes f_j^* = (\sum_{k_1} u_{k_1 i} e_{k_1}) \otimes (\sum_{k_2} v_{k_2 j} f_j)$$
$$= \sum_{k_1, k_2} u_{k_1 i} v_{k_2 j} (e_{k_1} \otimes f_{k_2}),$$

and

$$\lim_{i+j \to \infty} \|u_{k_1 i} v_{k_2 j}\| = 0,$$

we deduce that the row vectors of the transition matrix $U \otimes V$ are convergent sequences. Similarly, the row vectors of the inverse transition matrix $U^{-1} \otimes V^{-1}$ are convergent sequences. This proves that $\vec{e^*} \otimes \vec{f^*}$ is indeed a basis of the formal tensor product $M \otimes N$.

It remains to show that the nuclear map $\phi \otimes \psi$ on the formal tensor product $M \otimes N$ is independent of the choice of the basis. That is, we want to show that if

$$\sum_{i,j} a_{ij}(e_i \otimes f_j) = \sum_{i,j} a_{ij}^*(e_i^* \otimes f_j^*),$$

then we must have

$$\sum_{i,j} \sigma(a_{ij}) \phi(e_i) \otimes \psi(f_j) = \sum_{i,j} \sigma(a_{ij}^*) \phi(e_i^*) \otimes \psi(f_j^*).$$

Replacing $\vec{e^*}$ by $\vec{e}U$ (resp. $\vec{f^*}$ by $\vec{f}V$), we need to check that if

$$\sum_{i,j} a_{ij}(e_i \otimes f_j) = \sum_{k_1, k_2} (\sum_{i,j} a_{ij}^* u_{k_1 i} v_{k_2 j}) e_{k_1} \otimes f_{k_2},$$

then

$$\sum_{i,j} \sigma(a_{ij}) \phi(e_i) \otimes \psi(f_j) = \sum_{k_1, k_2} (\sum_{i,j} \sigma(a_{ij}^* u_{k_1 i} v_{k_2 j})) \phi(e_{k_1}) \otimes \psi(f_{k_2}).$$

But this follows from the infinite $\sigma$-linearity of $\phi$ and $\psi$. The lemma is proved. $\qquad \square$

COROLLARY 2.6.    *The category of nuclear $\sigma$-modules is closed under direct sums, (formal) tensor product, symmetric powers and exterior powers.*

The $\sigma$-linear map $\phi$ in a nuclear $\sigma$-module $(M, \phi)$ can also be viewed as an infinite $A_0$-linear map $\phi : M^{(\sigma)} \to M$ of Banach $A_0$-modules, where $M^{(\sigma)}$ is the (formal) pull back of $M$ by $\sigma : A_0 \to A_0$. That is,

$$M^{(\sigma)} = \{\sum_i b_i \otimes e_i \mid b_i \in A_0\},$$

where for $a_i \in A_0$, we have the relation

$$\sum \sigma(a_i) \otimes e_i = \sum_i a_i e_i.$$

The map $\phi$ is called the Frobenius map of the nuclear $\sigma$-module $M$. We shall often just write $M$ or $\phi$ for the pair $(M, \phi)$. To be more specific about the nuclear condition of the nuclear map $\phi$ with respect to a basis $\vec{e}$ of $M$, we shall introduce several notions attached to a given basis.

*Definition* 2.7. Let $(M, \phi)$ be a nuclear $\sigma$-module with a given basis $\vec{e}$. For each integer $1 \leq i < \infty$, let $d_i$ be the smallest positive integer $d$ such that for all $j > d$,

$$\phi(e_j) \equiv 0 \pmod{\pi^i}.$$

This is a finite integer for each $i$ since

$$\lim_j \|\phi(e_j)\| = 0.$$

For $0 \leq i < \infty$, define

$$h_i = d_{i+1} - d_i, \ d_0 = 0.$$

The sequence $h = h(\vec{e}) = \{h_0, h_1, \cdots\}$ is called the *basis sequence* of $\phi$ with respect to the basis $\vec{e}$. If we let $M_{(i)}$ be the Banach $A_0$-module with the basis $\{e_{d_i+1}, e_{d_i+2}, \cdots\}$, then there is a decreasing *basis filtration* of Banach $A_0$-submodules:

$$M = M_{(0)} \supset M_{(1)} \supset \cdots \supset M_{(j)} \supset \cdots \supset 0, \ \cap_j M_{(j)} = 0,$$

where each $M_{(i)}$ has a basis, each quotient $M_{(i)}/M_{(i+1)}$ is a finite free $A_0$-module and

$$\phi(M_{(i)}) \subseteq \pi^i M, \quad \mathrm{rank}(M_{(i)}/M_{(i+1)}) = h_i.$$

Equivalently, the matrix $B$ of $\phi$ with respect to $\vec{e}$ (defined by $\phi(\vec{e}) = \vec{e}B$) is of the form

$$B = (B_0, \pi B_1, \pi^2 B_2, \cdots, \pi^i B_i, \cdots),$$

where each block $B_i$ is a matrix over $A_0$ with $h_i$ (finitely many) columns.

We now turn to discussing the convergent condition of $\phi$ as a "power series" in term of $X$, generalizing the overconvergent and $c \log$-convergent condition of an element in $A_0$.

*Definition* 2.8. Let $B(X)$ be the matrix of a nuclear $\sigma$-module $(M, \phi)$ under a basis $\vec{e}$. Write

$$B(X) = \sum_{u \in \mathbf{Z}^n} B_u X^u,$$

where the coefficients $B_u$ are (possibly infinite) matrices with entries in $R$. We say that $(M, \phi)$ is *convergent* if the matrix $B(X)$ is convergent. Namely, $B(X)$ satisfies the condition

$$(2.3) \qquad \lim_{|u| \to \infty} \operatorname{ord}_\pi B_u = \infty.$$

Similarly, we say that $(M, \phi)$ is *overconvergent* if the matrix $B(X)$ is overconvergent. Namely, $B(X)$ satisfies the condition

$$(2.4) \qquad \lim_{|u| \to \infty} \inf \frac{\operatorname{ord}_\pi B_u}{|u|} > 0.$$

Likewise, we say that $(M, \phi)$ is $c$ log-*convergent* if the matrix $B(X)$ is $c$ log-convergent. Namely, $B(X)$ satisfies the condition

$$(2.5) \qquad \lim_{|u| \to \infty} \inf \frac{\operatorname{ord}_\pi B_u}{\log_q |u|} \geq c.$$

Note that our definition of the convergence property depends on the choice of the basis $\vec{e}$. It may happen that the matrix $B(X)$ under one basis is overconvergent (resp. $c$ log-convergent) but the matrix under a new basis is *not* overconvergent (resp. $c$ log-convergent). Our definition is that as long as there is one basis for which the matrix is overconvergent (resp. $c$ log-convergent), then the $\sigma$-module $(M, \phi)$ is called overconvergent (resp., $c$ log-convergent). In the case of finite rank, an overconvergent $\sigma$-module is then a $\sigma$-module which can be defined over $A^\dagger$ and then by base extension $A^\dagger \to A_0$. In the case of infinite rank, an overconvergent nuclear $\sigma$-module is stronger than a nuclear $\sigma$-module which can be defined over $A^\dagger$. It may happen that the matrix $B(X) = (b_{i,j}(X))$ has all of its entries $b_{i,j}(X)$ in $A^\dagger$ without satisfying the uniformity requirement (2.4). Similar statements hold for $c$ log-convergent nuclear $\sigma$-modules as well.

The matrix $B(X)$ of a nuclear $\sigma$-module under a basis can be written in the form

$$B(X) = \sum_{u \in \mathbf{Z}^n} B_u X^u, \quad B_u = (b_{w_1, w_2}(u)),$$

where $w_1$ denotes the row index and $w_2$ denotes the column index of the constant matrix $B_u$. Conversely, such a matrix power series $B(X)$ is the matrix of some $c$ log-convergent nuclear $\sigma$-module under some basis if and only if $B(X)$ satisfies the following two conditions. First, we have the $c$ log-convergent condition for $B(X)$ as a power series in $X$:

$$\lim_{|u| \to \infty} \inf \frac{\operatorname{ord}_\pi B_u}{\log_q |u|} \geq c.$$

Second, we have the nuclear condition for $B(X)$ to be the matrix of $\phi$ under a nuclear basis. This condition can be rephrased as follow: For any positive number $C > 0$, there is an integer $N_C > 0$ such that for all column indexes $w_2 > N_C$, we have

$$\mathrm{ord}_\pi b_{w_1,w_2}(u) \geq C,$$

uniformly for all $u$ and all $w_1$.

As stated in Corollary 2.6, the category of nuclear $\sigma$-modules is closed under direct sum, (formal) tensor product, symmetric product and exterior product. The same statement holds for the category of convergent nuclear $\sigma$-modules (resp. $c$ log-convergent nuclear $\sigma$-modules, resp. overconvergent nuclear $\sigma$-modules). The category of convergent nuclear $\sigma$-modules includes all $\sigma$-modules of finite rank. In the case that the finite rank map $\phi : M^{(\sigma)} \to M$ is an isogeny, namely when $\phi : M^{(\sigma)} \otimes \mathbf{Q} \to M \otimes \mathbf{Q}$ is an isomorphism of $A_0 \otimes \mathbf{Q}$-modules, the $\sigma$-module $(M, \phi)$ more or less corresponds to an F-crystal on the toric space $\mathbf{G}_m^n$. If, in addition, the finite rank map $\phi : M^{(\sigma)} \to M$ is an isomorphism of $A_0$-modules, the $\sigma$-module $(M, \phi)$ corresponds exactly to a unit root F-crystal on $\mathbf{G}_m^n$. The implied condition giving the horizontal connection can be easily provided. As Nicholas Katz puts it, the connection comes for free. The notion of overconvergent nuclear $\sigma$-modules of finite rank arises naturally from various relative $p$-adic (but not $p$-adic étale) cohomology groups of a family of varieties over finite fields. The notion of overconvergent nuclear $\sigma$-modules of infinite rank already arises if one works on the chain level of $p$-adic cohomology groups of a family of varieties over a finite field.

From now on, we shall exclusively work with a nuclear $\sigma$-module $(M, \phi)$ and often identify it with its matrix $B(X)$ under a basis. We can define an $L$-function for any nuclear $\sigma$-module. For a geometric point $\bar{x} \in \mathbf{G}_m^n$ of degree $d$ over $\mathbf{F}_q$, let $x$ be the Teichmüller lifting of $\bar{x}$. The fibre $M_x$ at $x$ is a Banach module with a countable basis over $R[x]$, where $R[x]$ is the ring obtained from $R$ by adjoining the coordinates of $x$. Via the $R$-algebra homomorphism $R[X] \to R[x]$ defined by $X \to x$, the map $\sigma$ acts on $R[x]$ by $\sigma(x) = x^q$. The fibre map

$$\phi_x : M_x \longrightarrow M_x$$

is $\sigma$-linear (not linear in general). But its $d^{\text{th}}$ iterate $\phi_x^d$ is $R[x]$-linear (since $\sigma^d(x) = x$) and can be defined over $R$. Furthermore, $\phi_x^d$ is the dual of a completely continuous operator in the sense of Serre [20]. This implies that $\phi_x^d$ itself is nuclear over $R[x]$. Thus, the characteristic series $\det(I - \phi_x^d T)$ (or the Fredholm determinant) of $\phi_x^d$ acting on $M_x$ (descended from $R[x]$ to $R$) is well-defined with coefficients in $R$. It is a $p$-adic entire function and is independent of the choice of the geometric point $\bar{x}$ associated to the closed point over $\mathbf{F}_q$ containing $\bar{x}$. Thus, it makes sense to say "the characteristic power series of $\phi$ at a closed point".

*Definition* 2.9. The $L$-function of a nuclear $\sigma$-module $(M, \phi)$ is defined by the Euler product

$$L(\phi, T) = \prod_{\bar{x} \in \mathbf{G}_m^n / \mathbf{F}_q} \frac{1}{\det(I - \phi_x^{\deg(\bar{x})} T^{\deg(\bar{x})})},$$

where $\bar{x}$ runs over the closed points and $x$ is the Teichmüller lifting of $\bar{x}$.

Note that each of our Euler factors is the reciprocal of an entire power series, not necessarily a polynomial. We are working in a more general situation than that found in the existing literature. If $B(X)$ is the matrix of $\phi$ under a basis, then one checks that

$$\begin{aligned} L(\phi, T) &= L(B(X), T) \\ &= \prod_{\bar{x} \in \mathbf{G}_m^n / \mathbf{F}_q} \frac{1}{\det(I - B(x)B(x^q) \cdots B(x^{q^{\deg(\bar{x})-1}}) T^{\deg(\bar{x})})}. \end{aligned}$$

This shows that the $L$-function is a power series with coefficients in $R$. Thus, it is trivially meromorphic in the open unit disc $|T|_\pi < 1$. If $\phi$ is convergent (which is always the case if $M$ is of finite rank), it can be shown that the $L$-function is meromorphic on the closed unit disc $|T|_\pi \leq 1$. Dwork's fundamental trace formula can be extended to show that this $L$-function is $p$-adic meromorphic if $\phi$ is overconvergent. The optimal result is given by the following theorem, which was proved in [10] and [26] for the finite rank case.

THEOREM 2.10. *Let $(M, \phi)$ be a $c$ log-convergent nuclear $\sigma$-module for some $0 < c \leq \infty$. Then the $L$-function $L(\phi, T)$ is $p$-adic meromorphic in the open disc $|T|_\pi < p^c$, where $|\pi|_\pi = 1/p$.*

Actually, what we will need in this paper is a stronger family version of the above result as given later in Section 5. Dwork's conjecture says that for an invertible $(M, \phi)$ of finite rank embedded (see §3) in an ordinary overconvergent $\sigma$-module of finite rank, the $L$-function $L(\phi, T)$ is meromorphic. No single nontrivial example was known previously, other than the "trivial" $\infty$ log-convergent case already covered by the above theorem.

In order to describe Dwork's conjecture fully, we need to define a slightly more general $L$-function, which could be called a higher power zeta function. If $\phi$ is a nuclear $\sigma$-module and $k$ is a positive integer, we define the $k^{\text{th}}$ power $\phi^k$ to be $\phi$ iterated $k$-times. The pair $(M, \phi^k)$ can be viewed as a $\sigma^k$-module but it is *not* a $\sigma$-module as $\phi^k$ is $\sigma^k$-linear, not $\sigma$-linear, in general. The $L$-function of the $k^{\text{th}}$ power $\phi^k$ is defined to be

$$L(\phi^k, T) = \prod_{\bar{x} \in \mathbf{G}_m^n / \mathbf{F}_q} \frac{1}{\det(I - \phi_x^{k\deg(\bar{x})} T^{\deg(\bar{x})})}.$$

If $B(X)$ is the matrix of $\phi$ under a basis, then one checks that

$$L(\phi^k, T) = \prod_{\bar{x} \in \mathbf{G}_m^n / \mathbf{F}_q} \frac{1}{\det(I - B(x)B(x^q) \cdots B(x^{q^{k \deg(\bar{x})-1}}) T^{\deg(\bar{x})})}$$

$$= \prod_{\bar{x} \in \mathbf{G}_m^n / \mathbf{F}_q} \frac{1}{\det(I - (B(x)B(x^q) \cdots B(x^{q^{\deg(\bar{x})-1}}))^k T^{\deg(\bar{x})})}.$$

In other words, $L(\phi^k, T)$ is obtained from $L(\phi, T)$ by raising each recipro-cal root $\alpha$ of the Euler factors to its $k^{\mathrm{th}}$ power $\alpha^k$. We want to emphasize that $L(\phi^k, T)$ is different from the following function $L(B^k, T)$ which is the $L$-function of the nuclear $\sigma$-module with matrix $B^k$:

$$L(B^k, T) = \prod_{\bar{x} \in \mathbf{G}_m^n / \mathbf{F}_q} \frac{1}{\det(I - B^k(x)B^k(x^q) \cdots B^k(x^{q^{\deg(\bar{x})-1}}) T^{\deg(\bar{x})})}.$$

They are different because matrix multiplication is in general noncommutative if $B$ has rank greater than one.

Let $(M_1, \phi_1)$ and $(M_2, \phi_2)$ be two nuclear $\sigma$-modules. Let $k_1$ and $k_2$ be two positive integers. Although their powers $\phi_1^{k_1}$ and $\phi_2^{k_2}$ do not have a $\sigma$-module structure in general, we can define the $L$-function of the "tensor product" of these powers as follows:

$$L(\phi_1^{k_1} \otimes \phi_2^{k_2}, T) = \prod_{\bar{x} \in \mathbf{G}_m^n / \mathbf{F}_q} \frac{1}{\det(I - \phi_{1x}^{k_1 \deg(\bar{x})} \otimes \phi_{2x}^{k_2 \deg(\bar{x})} T^{\deg(\bar{x})})}.$$

That is, the Euler factor at $x$ of $L(\phi_1^{k_1} \otimes \phi_2^{k_2}, T)$ is the "tensor product" of the corresponding Euler factors at $x$ of $L(\phi_1^{k_1}, T)$ and $L(\phi_2^{k_2}, T)$. The above definition works for the negative integer $k$ as well, provided that $\phi$ is injective fibre by fibre and of finite rank. If $(M, \phi)$ is of rank one, then $L(\phi^k, T)$ is just the $L$-function $L(\phi^{\otimes k}, T)$ of the $k^{\mathrm{th}}$ tensor (or symmetric) power of $(M, \phi)$. If $(M, \phi)$ is of rank greater than one, then the $k^{\mathrm{th}}$ power $\phi^k$ $(k > 1)$ does not have a $\sigma$-module interpretation as matrix multiplication is noncommutative in general. Nevertheless, by the Newton-Waring formula, we see that $\phi^k$ can be viewed as a virtual $\sigma$-module, namely, an integral combination of $\sigma$-modules with both positive and negative coefficients. This implies that $L(\phi^k, T)$ can be written as an alternating product of $L$-functions of various $\sigma$-modules. We shall need something finer than the Newton-Waring formula in order to prove uniform results.

## 3. Hodge-Newton decomposition and Dwork's conjecture

Let $(M, \phi)$ be a nuclear $\sigma$-module with a basis $\vec{e}$. First, we briefly recall the construction of formal exterior powers and formal symmetric powers. Let $1 \leq i < \infty$. The $i^{\text{th}}$ (formal) exterior power $(\wedge^i M, \wedge^i \phi)$ is the nuclear $\sigma$-module with the basis

$$\wedge^i \vec{e} = \{\cdots, e_{k_1} \wedge \cdots \wedge e_{k_i}, \cdots\}, \ k_1 < k_2 < \cdots < k_i$$

and with the infinite $\sigma$-linear endomorphism $\wedge^i \phi$ defined by

$$\wedge^i(\phi)(e_{k_1} \wedge \cdots \wedge e_{k_i}) = \phi(e_{k_1}) \wedge \cdots \wedge \phi(e_{k_i}).$$

For $i = 0$, we define $(\wedge^0 M, \wedge^0 \phi)$ to be the trivial rank one $\sigma$-module $(A_0, \sigma)$. Similarly, the $i^{\text{th}}$ (formal) symmetric power $(\text{Sym}^i M, \text{Sym}^i \phi)$ is the nuclear $\sigma$-module with the basis

$$\text{Sym}^i \vec{e} = \{\cdots, e_{k_1} \cdots e_{k_i}, \cdots\}, \ k_1 \leq k_2 \leq \cdots \leq k_i$$

and with the infinite $\sigma$-linear endomorphism $\text{Sym}^i \phi$ defined by

$$\text{Sym}^i(\phi)(e_{k_1} \cdots e_{k_i}) = \phi(e_{k_1}) \cdots \phi(e_{k_i}).$$

The (formal) tensor power $(M^{\otimes i}, \phi^{\otimes i})$ is similarly constructed. We shall use the convention that $\text{Sym}^k \phi = 0$ for all negative integers $k < 0$.

Next, we introduce the basis polygon and Newton polygon of a nuclear $\sigma$-module $(M, \phi)$ with a basis $\vec{e}$. Denote by $\text{ord}(\phi)$ to be the greatest integer $b$ such that

$$\phi \equiv 0 \ (\text{mod} \pi^b).$$

We say that $(M, \phi)$ is divisible by $\pi$ if $\phi$ is divisible by $\pi$. If $\phi$ is divisible by $\pi$, then we can write $\phi = \pi \phi_1$. In this case, one checks that

$$L(\phi, T) = L(\phi_1, \pi T).$$

Thus, for our purpose of studying $L$-functions, we may assume that $(M, \phi)$ is not divisible by $\pi$.

*Definition* 3.1. Let $(M, \phi)$ be a nuclear $\sigma$-module with a basis $\vec{e}$. Let $h = h(\vec{e}) = (h_0, h_1, \cdots)$ be the basis sequence of $\phi$ with respect to the basis $\vec{e}$ (see Definition 2.7). We define the *basis polygon* $P(\vec{e})$ of $(M, \phi)$ with respect to $\vec{e}$ to be the convex closure in the plane of the following lattice points:

$$(0,0), (h_0, 0), (h_0 + h_1, h_1), (h_0 + h_1 + h_2, h_1 + 2h_2), \cdots,$$

$$(h_0 + \cdots + h_i, h_1 + 2h_2 + \cdots + ih_i), \cdots.$$

That is, the basis polygon $P(\vec{e})$ is the polygon with a side of slope $i$ and horizontal length $h_i$ for every integer $0 \leq i < \infty$.

Note that our definition of the basis polygon $P(\vec{e})$ and the basis sequence $h(\vec{e})$ depends on the given basis $\vec{e}$. The basis polygon and the basis sequence are somewhat similar to the Hodge polygon and Hodge numbers. But they are different in general, even in the finite rank case. The same $\sigma$-module with different bases will give rise to different basis polygons and different basis sequences. The (generic) abstract Hodge polygon in the finite rank case as discussed in [15] does not depend on the given filtration. One can show that the generic Hodge polygon lies on or above our basis polygon. In the ordinary case to be defined below, the basis $\vec{e}$ will be the optimal one and our basis polygon is the same as the Hodge polygon in the finite rank case. Since we are dealing with the ordinary case of Dwork's conjecture in this paper, nothing is lost. It may be interesting to define an optimal intrinsic (necessarily generic) basis polygon in general which coincides with the Hodge polygon in the finite rank case. We shall not pursue it here. Our simplified definition here is just for the purpose of quickly getting to Dwork's conjecture without losing the key property we need to use.

*Definition* 3.2. For each closed point $\bar{x} \in \mathbf{G}_m^n$ over $\mathbf{F}_q$, the Newton polygon of the nuclear $\sigma$-module $(M, \phi)$ at $\bar{x}$ is the Newton polygon of the entire characteristic series $\det(I - \phi_x^{\deg(\bar{x})}T)$ defined with respect to the valuation $\mathrm{ord}_{\pi^{\deg(\bar{x})}}$.

A standard argument shows that the Newton polygon at each $\bar{x}$ lies on or above any basis polygon $P(\vec{e})$. This is known to be true in the finite rank case. Since we do not need it in this paper, we do not include the detail of the proof here

*Definition* 3.3. Let $\vec{e}$ be a basis of a nuclear $\sigma$-module $(M, \phi)$. The basis $\vec{e}$ is called *ordinary* if the basis polygon $P(\vec{e})$ of $(M, \phi)$ coincides with the Newton polygon of each fibre $(M, \phi)_x$, where $x$ is the Teichmüller lifting of the closed point $\bar{x} \in \mathbf{G}_m^n / \mathbf{F}_q$. The basis filtration (see Definition 2.7) attached to $\vec{e}$ is called ordinary if $\vec{e}$ is ordinary. Similarly, the basis $\vec{e}$ is called ordinary up to slope $j$ if the basis polygon $P(\vec{e})$ coincides with the Newton polygon of each fibre for all sides up to slope $j$. In particular, $\vec{e}$ is said to be ordinary at slope zero if the horizontal side of the Newton polygon at each fibre is of length $h_0(\vec{e})$. We say that $(M, \phi)$ is ordinary (resp. ordinary up to slope $j$, resp. ordinary at slope zero) if it has a basis which is ordinary (resp. ordinary up to slope $j$, resp. ordinary at slope zero).

An important property is that the category of ordinary nuclear $\sigma$-modules is closed under direct sum, tensor product, symmetric power and exterior power. This property is not hard to prove. It follows from the proof of the Hodge-Newton decomposition to be described below. Alternatively, one could check directly. For example, let $(M_1, \phi_1)$ be a nuclear $\sigma$-module with a basis $\vec{e}$ and with its associated basis sequence $h = (h_0, h_1, \cdots)$. Let $(M_2, \phi_2)$ be a nuclear $\sigma$-module with a basis $\vec{f}$ and with its associated basis sequence $g = (g_0, g_1, \cdots)$. Then the direct sum $(M_1 \oplus M_2, \phi_1 \oplus \phi_2)$ is a nuclear $\sigma$-module with the basis $\vec{e} \oplus \vec{f}$. Its associated basis sequence is given by

$$h \oplus g = (h_0 + g_0, h_1 + g_1, \cdots).$$

The tensor product $(M_1 \otimes M_2, \phi_1 \otimes \phi_2)$ is a nuclear $\sigma$-module with basis $\vec{e} \otimes \vec{f}$. Its associated basis sequence is given by

$$h \otimes g = (h_0 g_0, h_0 g_1 + h_1 g_0, h_0 g_2 + h_1 g_1 + h_2 g_0, \cdots).$$

The symmetric square $(\mathrm{Sym}^2 M_1, \mathrm{Sym}^2 \phi_1)$ is a nuclear $\sigma$-module with basis $\mathrm{Sym}^2 \vec{e}$. Its associated basis sequence $\mathrm{Sym}^2 h$ is given by

$$\left( \frac{h_0^2 + h_0}{2}, h_0 h_1, h_0 h_2 + \frac{h_1^2 + h_1}{2}, h_0 h_3 + h_1 h_2, h_0 h_4 + h_1 h_3 + \frac{h_2^2 + h_2}{2}, \cdots \right).$$

The exterior square $(\wedge^2 M_1, \wedge^2 \phi_1)$ is a nuclear $\sigma$-module with basis $\wedge^2 \vec{e}$. Its associated basis sequence is given by

$$\wedge^2 h = \left( \frac{h_0^2 - h_0}{2}, h_0 h_1, h_0 h_2 + \frac{h_1^2 - h_1}{2}, h_0 h_3 + h_1 h_2, h_0 h_4 + h_1 h_3 + \frac{h_2^2 - h_2}{2}, \cdots \right).$$

If $\vec{e}$ is an ordinary basis of $(M_1, \phi_1)$ and if $\vec{f}$ is an ordinary basis of $(M_2, \phi_2)$, one checks that $\vec{e} \oplus \vec{f}, \vec{e} \otimes \vec{f}, \mathrm{Sym}^2 \vec{e}$ and $\wedge^2 \vec{e}$ are, respectively, an ordinary basis of $(M_1 \oplus M_2, \phi_1 \oplus \phi_2)$, $(M_1 \otimes M_2, \phi_1 \otimes \phi_2)$, $(\mathrm{Sym}^2 M_1, \mathrm{Sym}^2 \phi_1)$ and $(\wedge^2 M_1, \wedge^2 \phi_1)$. The associated basis sequence is respectively $h \oplus g$, $h \otimes g$, $\mathrm{Sym}^2 h$ and $\wedge^2 h$. A systematic study of Newton polygons and Hodge polygons in finite rank case is given in [15].

Since $A_0$ is complete, any contraction mapping from $A_0$ into itself has a fixed point. This fact is all one needs to prove the following extension of Hodge-Newton decomposition, first proved by Dwork [8] in the finite rank case. The infinite rank case had already been considered and used in [23] to study the Newton polygon of Fredholm determinants arising from $L$-functions of exponential sums.

LEMMA 3.4. *Let $(M, \phi)$ be a nuclear $\sigma$-module ordinary at slope zero. Then there is a free $A_0$-submodule $M_0$ of rank $h_0$ transversal to the ordinary basis filtration: $M = M_0 \oplus M_{(1)}$, such that $\phi$ is stable on $M_0$ and $(M_0, \phi)$ is a unit root $\sigma$-module of rank $h_0$.*

*Proof.* The proof is completely similar to Dwork's original proof in the finite rank case. We give a brief sketch of the proof in order to have a feeling about what is going on. Under a basis $\vec{e}$ ordinary at slope zero, the matrix $B$ of $\phi$ is of the form

$$B = \begin{pmatrix} B_{00} & \pi B_{01} \\ B_{10} & \pi B_{11} \end{pmatrix},$$

where $B_{00}$ is an $h_0 \times h_0$ matrix over $A_0$. Since the basis $\vec{e}$ of $(M, \phi)$ is ordinary at slope zero, the matrix $B_{00}$ is invertible over $A_0$. Consider the new basis

$$(e_1, \cdots, e_i, \cdots) \begin{pmatrix} I_{00} & 0 \\ E_{10} & I_{11} \end{pmatrix} = \vec{e}E$$

of $M$, where $I_{00}$ is the identity matrix of rank $h_0$, $I_{11}$ is an identity matrix and $E_{10}$ is some matrix to be determined. The matrix $E_{10}$ has $h_0$ columns. One calculates that the matrix of $\phi$ under the new basis $\vec{e}E$ is given by

$$E^{-1}BE^{\sigma} = \begin{pmatrix} B_{00} + \pi B_{01}E_{10}^{\sigma} & \pi B_{01} \\ -E_{10}B_{00} + B_{10} + \pi(B_{11} - E_{10}B_{01})E_{10}^{\sigma} & \pi(B_{11} - E_{10}B_{01}) \end{pmatrix}.$$

To prove the theorem, we need to choose $E_{10}$ suitably such that the above new matrix is triangular. Namely, we want

$$-E_{10}B_{00} + B_{10} + \pi(B_{11} - E_{10}B_{01})E_{10}^{\sigma} = 0.$$

This equation can be rewritten as

$$E_{10} = B_{10}B_{00}^{-1} + \pi(B_{11} - E_{10}B_{01})E_{10}^{\sigma}B_{00}^{-1}.$$

It defines a contraction map. Thus, by successive iteration or the fixed point theorem in a $p$-adic Banach space, we conclude that there is a unique solution matrix $E_{10}$ with entries in the ring $A_0$. Furthermore, one checks that the map in the above change of basis satisfies the condition of Lemma 2.1. Thus, $\vec{e}E$ is indeed a basis with the same $h_0$. The lemma is proved. We note that the above change of basis will in general destroy any overconvergent property the original matrix $B(X)$ might have because the solution matrix $E_{10}$ will only be convergent even though $B$ is overconvergent.

In the situation of Lemma 3.4, we shall denote the sub $\sigma$-module $(M_0, \phi)$ by $(U_0, \phi_0)$. We say that the invertible $\sigma$-module $(U_0, \phi_0)$ is embedded in the ambient ordinary $\sigma$-module $(M, \phi)$. Alternatively, we shall say that $(U_0, \phi_0)$ is the unit root part (or the invertible part) of $(M, \phi)$. In the geometric situation, the unit root part can be constructed using expansion coefficients of differential forms [16]. It can also be obtained as the relative $p$-adic ètale cohomology with compact support of a family of varieties defined over a finite field. Dwork's unit root conjecture, extended to our setting, is then the following:

UNIT ROOT CONJECTURE.    *Let $(M, \phi)$ be an overconvergent nuclear $\sigma$-module ordinary at slope zero. Let $(U_0, \phi_0)$ be the unit root part of $(M, \phi)$. Then for every integer $k$, the unit root zeta function $L(\phi_0^k, T)$ is $p$-adic meromorphic.*

*Remark.* Even though the ambient $\sigma$-module $(M, \phi)$ is assumed to be overconvergent, its unit root part (obtained by solving the fixed point of a contraction map) will no longer be overconvergent in general. In fact, one can only expect the unit root part to be $c$ log-convergent for some small real number $c$. Dwork and Sperber showed in [10] that one can take $c=(p-1)/(p+1)$ in the finite rank case. They suggested the possibility that one might be able to take $c = 1$. This was later proved to be true by Sperber-Wan (unpublished). It is not very hard to construct examples to show that $c = 1$ is best possible.

A weaker version of our main result is the following:

UNIT ROOT THEOREM. *Let $(M, \phi)$ be a $c$ log-convergent nuclear $\sigma$-module $(M, \phi)$ for some $0 < c \leq \infty$, ordinary at slope zero. Let $(U_0, \phi_0)$ be the unit root part of $(M, \phi)$. Assume that $\phi_0$ has rank one. Then for each integer $k$, the unit root zeta function $L(\phi_0^k, T)$ is $p$-adic meromorphic in the open disc $|T|_\pi < p^c$.*

COROLLARY 3.5. *Let $(M, \phi)$ be an $\infty$ log-convergent nuclear $\sigma$-module $(M, \phi)$ ordinary at slope zero. Let $(U_0, \phi_0)$ be the unit root part of $(M, \phi)$. Assume that $\phi_0$ has rank one. Then for each integer $k$, the unit root zeta function $L(\phi_0^k, T)$ is $p$-adic meromorphic everywhere.*

This corollary shows that Dwork's unit root conjecture is true if the unit root part has rank one. Note that our result is actually stronger, since we only assume the ambient $\sigma$-module to be $c$ log convergent instead of being overconvergent. The $c$ log-convergence condition in the above theorem is best possible by the counterexample in [26]. The possibility of replacing the overconvergent condition by the weaker $\infty$ log-convergent condition in Dwork's conjecture (if true) was already suggested in [25]. All the above results can be extended to higher slope. The following Hodge-Newton decomposition is also due to Dwork [8] in the finite rank case. It can be proved in a similar way and uses induction.

LEMMA 3.6. *If a nuclear $\sigma$-module $(M, \phi)$ is ordinary up to slope $j$ for some integer $j \geq 0$, then there is an increasing $\phi$-stable filtration of free $A_0$-submodules of finite rank,*

$$0 \subset M_0 \subset M_1 \subset \cdots \subset M_j \subset M$$

*which is transversal to the ordinary basis filtration: $M = M_i \oplus M_{(i+1)}$ for all $0 \leq i \leq j$. Furthermore, for $0 \leq i \leq j$, the quotient $(M_i/M_{i-1}, \phi)$ is a $\sigma$-module of the form $(U_i, \pi^i \phi_i)$, where $(U_i, \phi_i)$ is a unit root $\sigma$-module of rank $h_i$.*

We shall call $(U_i, \phi_i)$ the unit root $\sigma$-module coming from the slope $i$ part of $(M, \phi)$. The inductive proof shows that if $(M, \phi)$ is ordinary up to slope $i$, then

$$(3.1) \qquad \operatorname{ord}(\wedge^{h_0 + \cdots + h_i} \phi) = h_1 + 2h_2 + \cdots + ih_i.$$

Dwork [8] also made the following more general conjecture.

HIGHER SLOPE CONJECTURE.    *Let $(M, \phi)$ be an overconvergent nuclear $\sigma$-module $(M, \phi)$ ordinary up to slope $j$ for some integer $j \geq 0$. Let $(U_i, \phi_i)$ be the unit root $\sigma$-module coming from the slope $i$ part of $(M, \phi)$, where $0 \leq i \leq j$. Then for each $0 \leq i \leq j$ and each integer $k$, the unit root zeta function $L(\phi_i^k, T)$ is $p$-adic meromorphic.*

Similarly, we can prove the following result for a higher slope zeta function.

HIGHER SLOPE THEOREM.    *Let $(M, \phi)$ be a $c$ log-convergent nuclear $\sigma$-module $(M, \phi)$ for some $0 < c \leq \infty$, ordinary up to slope $j$ for some integer $j \geq 0$. Let $(U_j, \phi_j)$ be the unit root $\sigma$-module coming from the slope $j$ part of $(M, \phi)$. Assume that $\phi_j$ has rank one. Then for each integer $k$, the unit root zeta function $L(\phi_j^k, T)$ is $p$-adic meromorphic in the open disc $|T|_\pi < p^c$.*

COROLLARY 3.7.    *Let $(M, \phi)$ be an $\infty$ log-convergent nuclear $\sigma$-module $(M, \phi)$ ordinary up to slope $j$ for some integer $j \geq 0$. Let $(U_j, \phi_j)$ be the unit root $\sigma$-module coming from the slope $j$ part of $(M, \phi)$. Assume that $\phi_j$ has rank one. Then for each integer $k$, the unit root zeta function $L(\phi_j^k, T)$ is $p$-adic meromorphic everywhere.*

## 4. *L-functions of various tensor powers*

Let $E$ be a finite square matrix over a field of characteristic zero. Our first goal is to prove a good universal formula which expresses the characteristic polynomial $\det(I - TE^k)$ of the $k^{\text{th}}$ power $E^k$ in terms of a certain alternating product of the characteristic polynomials of various tensor and exterior products of $E$. There are several such formulas available. The most well-known one is the classical Newton-Waring formula expressing the power symmetric functions in terms of elementary symmetric functions. To get uniform results, we will need something finer than the Newton-Waring formula. Since our formulas are actually universal formulas, they hold for fields of positive characteristic as well. For simplicity of our proof, we shall assume that we are in characteristic zero. Recall that $\mathrm{Tr}(E)$ denotes the trace of the matrix $E$. We make the convention that $\mathrm{Sym}^k E$ is the zero matrix (of any size) for negative integer $k < 0$.

LEMMA 4.1.    *Let $E$ be a finite square matrix. Then, for each integer $k > 0$, there exists the relation*

$$(4.1) \qquad \mathrm{Tr}(E^k) = \sum_{i \geq 1} (-1)^{i-1} i \ \mathrm{Tr}(\mathrm{Sym}^{k-i} E) \ \mathrm{Tr}(\wedge^i E).$$

*Proof.* This can be proved using a combinatorial counting argument. We shall, however, use the following more conceptual proof communicated to us by N. Katz. Taking the logarithmic derivative of the well-known identity

$$\frac{1}{\det(I - ET)} = \exp(\sum_{k \geq 1} \frac{\mathrm{Tr}(E^k)}{k} T^k),$$

we get

$$\frac{-T\frac{d}{dT}\det(I - ET)}{\det(I - ET)} = \sum_{k \geq 1} \mathrm{Tr}(E^k)T^k.$$

On the other hand, by looking at the eigenvalues, one easily gets

$$(4.2) \qquad \det(I - ET) = \sum_{i \geq 0} (-1)^i \, \mathrm{Tr}(\wedge^i E)T^i,$$

$$(4.3) \qquad \frac{1}{\det(I - ET)} = \sum_{j \geq 0} \mathrm{Tr}(\mathrm{Sym}^j E)T^j.$$

It follows that

$$(\sum_{i \geq 0} (-1)^{i-1} i \, \mathrm{Tr}(\wedge^i E)T^i)(\sum_{j \geq 0} \mathrm{Tr}(\mathrm{Sym}^j E)T^j) = \sum_{k \geq 1} \mathrm{Tr}(E^k)T^k.$$

Comparing the coefficients of $T^k$ on both sides, we conclude that the lemma is true. $\qquad \square$

The next result is a further refinement in the sense that the coefficients have smaller size. This would be useful in practical calculation of the unit root zeta function.

LEMMA 4.2. *Let $E$ be a finite square matrix. Then, for each integer $k > 0$,*

$$\mathrm{Tr}(E^k) = \mathrm{Tr}(\mathrm{Sym}^k E) + \sum_{i \geq 2} (-1)^{i-1}(i - 1) \, \mathrm{Tr}(\mathrm{Sym}^{k-i} E) \, \mathrm{Tr}(\wedge^i E).$$

*Proof.* Using the trivial identity

$$\frac{\det(I - ET)}{\det(I - ET)} = 1$$

and equations (4.2)–(4.3), we get

$$(\sum_{i \geq 0} (-1)^i \, \mathrm{Tr}(\wedge^i E)T^i)(\sum_{j \geq 0} \mathrm{Tr}(\mathrm{Sym}^j E)T^j) = 1.$$

Comparing the coefficients of $T^k$, we deduce that for $k \geq 1$,

$$\sum_{i \geq 0} (-1)^i \, \mathrm{Tr}(\mathrm{Sym}^{k-i} E) \, \mathrm{Tr}(\wedge^i E) = 0.$$

Adding this equation to (4.1), we immediately get Lemma 4.2. $\qquad \square$

We now apply the above two lemmas to the tensor product $E_1^{k_1} \otimes E_2^{k_2}$. Since

$$\mathrm{Tr}(E_1^{k_1})\, \mathrm{Tr}(E_2^{k_2}) = \mathrm{Tr}(E_1^{k_1} \otimes E_2^{k_2}),$$

we deduce from Lemma 4.1 the following result.

LEMMA 4.3.    *Let $E_1$ and $E_2$ be two finite square matrices. Then, for integers $k_1 > 0$ and $k_2 > 0$,*

$$\mathrm{Tr}(E_1^{k_1} \otimes E_2^{k_2}) = \sum_{i_1, i_2 \geq 1} (-1)^{i_1 + i_2} i_1 i_2 \, \mathrm{Tr}(\mathrm{Sym}^{k_1 - i_1} E_1) \, \mathrm{Tr}(\wedge^{i_1} E_1)$$
$$\times \quad \mathrm{Tr}(\mathrm{Sym}^{k_2 - i_2} E_2) \, \mathrm{Tr}(\wedge^{i_2} E_2).$$

If one uses Lemma 4.2 instead of Lemma 4.1, one gets a more efficient but less compact formula for $\mathrm{Tr}(E_1^{k_1} \otimes E_2^{k_2})$. We omit it. The above results on matrix trace can now be transformed to similar results for characteristic polynomials.

LEMMA 4.4.    *Let $E$ be a finite square matrix. Then, for each integer $k > 0$,*

$$\det(I - TE^k) = \prod_{i \geq 1} \det(I - T(\mathrm{Sym}^{k-i} E \otimes \wedge^i E))^{(-1)^{i-1} i}.$$

*Proof.* For each $j \geq 1$, Lemma 4.1 shows that

$$\mathrm{Tr}(E^{jk}) = \sum_{i \geq 1} (-1)^{i-1} i \, \mathrm{Tr}(\mathrm{Sym}^{k-i}(E^j) \otimes \wedge^i(E^j))$$
$$= \sum_{i \geq 1} (-1)^{i-1} i \, \mathrm{Tr}((\mathrm{Sym}^{k-i} E \otimes \wedge^i E)^j).$$

These relations together with (4.3) immediately imply the lemma. The proof is complete.    □

More generally, we consider the tensor product of powers of two matrices. The following formula is a consequence of Lemma 4.3.

LEMMA 4.5.    *Let $E_1$ and $E_2$ be two finite square matrices. Then, for integers $k_1 > 0$ and $k_2 > 0$,*

$$\det(I - TE_1^{k_1} \otimes E_2^{k_2}) = \prod_{i_1, i_2} \det(I - TE(i_1, i_2))^{(-1)^{(i_1 + i_2)} i_1 i_2},$$

*where*

$$E(i_1, i_2) = \mathrm{Sym}^{k_1 - i_1} E_1 \otimes \wedge^{i_1} E_1 \otimes \mathrm{Sym}^{k_2 - i_2} E_2 \otimes \wedge^{i_2} E_2.$$

The above formulas hold for infinite matrices as well if every object involved is convergent. They are universal formulas; see [26]. In fact, the same proof works if we let the entries of $E$ be just indeterminates. Then, the above formulas hold as identities of two formal power series in those indeterminates and $T$. We now apply this to certain infinite matrices over $R$.

*Definition* 4.6. An infinite matrix $E$ over $R$ is called *nuclear* if its column vector $V_j$ goes to zero as the column index $j$ goes to infinity. Namely, $\lim_j \operatorname{ord}_\pi V_j = \infty$.

If $E$ is the transpose of the matrix of a completely continuous endomorphism of a Banach space over $R$ under a countable orthonormal row basis, then $E$ is nuclear. The converse is also true. This shows that the characteristic series $\det(I - TE)$ is well-defined and entire for a nuclear matrix $E$. The set of nuclear matrices is clearly stable under direct sum, tensor product, exterior product and symmetric product. Thus, we have the following specialization from the above mentioned universal formula. It can also be proved directly from Lemma 4.4 by taking the limit.

LEMMA 4.7. *Let $E$ be a nuclear matrix over $R$. Then, for each integer $k > 0$,*
$$\det(I - TE^k) = \prod_{i \geq 1} \det(I - T(\operatorname{Sym}^{k-i} E \otimes \wedge^i E))^{(-1)^{i-1} i}.$$

If one uses Lemma 4.2 instead, one gets the following:

LEMMA 4.8. *Let $E$ be a nuclear matrix over $R$. Then, for each integer $k > 0$,*
$$\det(I - TE^k) = \det(I - T\operatorname{Sym}^k E) \times \prod_{i \geq 2} \det(I - T(\operatorname{Sym}^{k-i} E \otimes \wedge^i E))^{(-1)^{i-1}(i-1)}.$$

Similarly, there are two formulas for $\det(I - TE_1^{k_1} \otimes E_2^{k_2})$, where $E_1$ and $E_2$ are nuclear matrices. We omit them. We now apply Lemma 4.7 to the $L$-function of a nuclear $\sigma$-module to obtain some preliminary analytic information. For each integer $k \geq 1$, recall that we defined
$$L(\phi^k, T) = \prod_{\bar{x} \in \mathbf{G}_m^n / \mathbf{F}_q} \frac{1}{\det(I - \phi_x^{k \deg(\bar{x})} T^{\deg(\bar{x})})}.$$

With the above preparation, we can now prove a simple meromorphy result for the $k^{\text{th}}$ power $L$-function $L(\phi^k, T)$, which is again best possible by the example in [26].

LEMMA 4.9. *Let $(M, \phi)$ be a $c$ log-convergent nuclear $\sigma$-module for some $0 < c \leq \infty$. Then for each integer $k > 0$, the $L$-function $L(\phi^k, T)$ is $p$-adic meromorphic in the open disc $|T|_\pi < p^c$.*

*Proof.* Applying Lemma 4.7 (which is a universal formula) and standard representation theory, we find the following basic decomposition of $L$-functions:

$$(4.4) \qquad L(\phi^k, T) = \prod_{i \geq 1} L(\mathrm{Sym}^{k-i}\phi \otimes \wedge^i\phi, T)^{(-1)^{i-1}i}.$$

One can also derive (4.4) by applying Lemma 4.7 directly to each Euler factor in the definition of $L(\phi^k, T)$ and using the following "notational" identities:

$$(4.5) \qquad (\mathrm{Sym}^k\phi_x)^{\deg(x)} = \mathrm{Sym}^k(\phi_x^{\deg(x)}),$$

$$(4.6) \qquad (\wedge^k\phi_x)^{\deg(x)} = \wedge^k(\phi_x^{\deg(x)}).$$

In fact, both sides of (4.5) are the linear map $\phi_x^{\deg(x)}$ acting on the space $\mathrm{Sym}^k M_x$. The left side of (4.5) means the $\deg(x)^{\mathrm{th}}$ iterate of the semi-linear map $\phi_x$ acting on $(\mathrm{Sym}^k M)_x$. The right side of (4.5) means the $k^{\mathrm{th}}$ symmetric power of the linear map $\phi_x^{\deg(x)}$ acting on $M_x$, which is (by the definition of symmetric power) the linear map $\phi_x^{\deg(x)}$ acting on the $k^{\mathrm{th}}$ symmetric power $\mathrm{Sym}^k M_x$. This proves (4.5). The proof of (4.6) is very similar.

Note that the product in (4.4) is a finite product since $\mathrm{Sym}^{k-i}\phi = 0$ for $i > k$. Now, for each integer $i$, the map $\mathrm{Sym}^{k-i}\phi \otimes \wedge^i\phi$ is just the Frobenius map of the nuclear $\sigma$-module $\mathrm{Sym}^{k-i}M \otimes \wedge^i M$, which is $c$ log-convergent. By Theorem 2.4, each of the $L$-factors on the right side is $p$-adic meromorphic in the open disc $|T|_\pi < p^c$. Thus, the $L$-function on the left is also meromorphic in the same disc. The lemma is proved. $\qquad \square$

This result already gives some information about Dwork's unit root zeta function.

THEOREM 4.10. *Let $(U_0, \phi_0)$ be the unit root part of a $c$ log-convergent nuclear $\sigma$-module $(M, \phi)$ ordinary at slope zero. Then, for each integer $k > 0$, the unit root zeta function $L(\phi_0^k, T)$ is meromorphic in the open disc $|T|_\pi < p^{\min(k,c)}$.*

*Proof.* By Hodge-Newton decomposition (Lemma 3.4), there is a short exact sequence of nuclear $\sigma$-modules:

$$0 \longrightarrow M_0 \longrightarrow M \longrightarrow M_1 \longrightarrow 0$$

such that $(M_1, \phi_1)$ is divisible by $\pi$ and $(M_0, \phi_0)$ is a unit root $\sigma$-module of rank $h_0$. Write $\phi_1 = \pi\psi$. Then, we have the following decomposition of $L$-functions:

$$L(\phi^k, T) = L(\phi_0^k, T)L(\psi^k, \pi^k T).$$

Because of the shifting factor $\pi^k$, we see that the second factor $L(\psi^k, \pi^k T)$ on the right is trivially meromorphic in the open disc $|T|_\pi < p^k$. The factor $L(\phi^k, T)$ on the left is meromorphic in the open disc $|T|_\pi < p^c$ by Lemma 4.9. Thus, the first factor $L(\phi_0^k, T)$ on the right is meromorphic in the open disc $|T|_\pi < p^{\min(k,c)}$. The proof is complete.  $\square$

If the first nonzero slope of $\phi$ is $s_1$ instead of 1, the same proof shows that one gets the meromorphic continuation to the larger disc $|T|_\pi < p^{\min(s_1 k, c)}$. Actually, it can be improved a little bit further to the disc $|T|_\pi < p^{\min(s_1(k+1), c)}$. We shall, however, not pursue this direction here.

COROLLARY 4.11.    *Let* $(U_0, \phi_0)$ *be the unit root part of an* $\infty$ *log-convergent nuclear* $\sigma$*-module* $(M, \phi)$ *ordinary at slope zero. Then, for each integer* $k > 0$*, the unit root zeta function* $L(\phi_0^k, T)$ *is meromorphic in the open disc* $|T|_\pi < p^k$.

This simple result is the starting point for our work on Dwork's conjecture. It shows that we get better and better meromorphic continuation for $L(\phi_0^k, T)$ as $k$ grows. To get meromorphic continuation of $L(\phi_0^k, T)$ for a fixed $k$, we can try to use some sort of easily derived limiting formula such as

$$L(\phi_0^k, T) = \lim_{m \to \infty} L(\phi_0^{k + (q^{m!} - 1) p^m}, T).$$

In order to pass the better meromorphic property to its limit, as Deligne pointed out, we have to bound uniformly the number of zeros and poles of each member of the sequence of $L$-functions in any fixed finite disc. Namely, we need some type of uniformity result about meromorphy of the sequence of functions $L(\phi_0^k, T)$ for all large integers $k$. The desired uniformity does not necessarily hold in general, but does so if the rank of $\phi_0$ is one. In the next section, we shall prove the uniformity in the normalized rank one case.

## 5. Continuous and uniform families

In this section, we discuss continuous and uniform families of meromorphic functions. Let $S$ be a subset of a complete metric space. Let $\bar{S}$ denote the closure of $S$. That is, $\bar{S}$ is the union of $S$ and its limit points. In our applications, the space $S$ will be either an infinite subset of the $p$-adic integers $\mathbf{Z}_p$ with the induced $p$-adic topology, or an infinite subset of the integers $\mathbf{Z}$ with the sequence topology, where the sequence topology of an infinite subset $S = \{k_i\} \subset \mathbf{Z}$ indexed by positive integers $i \in \mathbf{Z}_{>0}$ means that the distance function on $S$ is given by

$$d(k_i, k_j) = |\frac{1}{i} - \frac{1}{j}|.$$

In the latter case, $\infty$ is the unique limiting point not contained in $S$. The set $S$ will be our parameter space for which the parameter $k$ varies. We shall only consider a family of power series in $R[[T]]$ of the following form:

$$f(k, T) = \sum_{m \geq 0} f_m(k) T^m, \quad f_0(k) = 1, \quad f_m(k) \in R,$$

where each coefficient $f_m(k)$ is a function from $S$ to $R$. Note that the constant term is always 1. Furthermore, the ring $R$ is always complete and compact with respect to the $p$-adic topology of $R$. The family $f(k, T)$ is called *continuous* in $k$ if each coefficient $f_m(k)$ is a continuous function from $S$ to $R$. The family $f(k, T)$ is called *uniformly continuous* in $k$ if each coefficient $f_m(k)$ is a uniformly continuous function from $S$ to $R$. That is, for each fixed $m$ and each $\epsilon > 0$, there is a real number $\delta(\epsilon) > 0$ such that whenever $k_1$ and $k_2$ are two elements of $S$ satisfying

$$d(k_1, k_2) < \delta(\epsilon),$$

we have

$$|f_m(k_1) - f_m(k_2)| < \epsilon.$$

For a real number $c \geq 0$, the family of functions $f(k, T)$ is called *uniformly analytic* in the open disc $|T|_\pi < p^c$ if

$$(5.1) \qquad\qquad \lim_{m \to \infty} \inf \frac{\inf_k \operatorname{ord}_\pi f_m(k)}{m} \geq c.$$

The family $f(k, T)$ is called *uniformly meromorphic* in the open disc $|T|_\pi < p^c$ if $f(k, T)$ can be written as a quotient $f_1(k, T)/f_2(k, T)$ of two families, where both $f_1(k, T)$ and $f_2(k, T)$ are uniformly analytic in the open disc $|T|_\pi < p^c$. It is clear that the product and quotient (if the denominator is not the zero family) of two uniformly meromorphic families (parametrized by the same parameter $k$) are still uniformly meromorphic.

If the set $S$ consists of a sequence of elements $\{k_i\}$, we shall call the family $f(k_i, T)$ a sequence of functions. The following result is an immediate consequence of the above defining inequality in (5.1).

LEMMA 5.1.    *Let* $f(k_i, T)$ *be a continuous sequence of power series. Assume that the sequence* $f(k_i, T)$ *is uniformly analytic in the open disc* $|T|_\pi < p^c$. *If* $k = \lim_i k_i$ *exists in* $\bar{S}$, *then the limit function*

$$g(k, T) = \lim_i f(k_i, T)$$

*exists and is analytic in the open disc* $|T|_\pi < p^c$. *Furthermore, if* $k \in S$, $g(k, T) = f(k, T)$.

COROLLARY 5.2.    *If* $f(k, T)$ *is a continuous family of uniformly analytic functions in* $|T|_\pi < p^c$ *for* $k \in S$, *then* $f(k, T)$ *extends uniquely to a family of* *uniformly analytic functions in $|T|_\pi < p^c$ for $k \in \bar{S}$ (the closure of $S$). If, in addition, the original family $f(k, T)$ parametrized by $k \in S$ is uniformly continuous, then the extended family $f(k, T)$ parametrized by $k \in \bar{S}$ is uniformly continuous.*

*Proof.* The first part is clear. The second part can be proved in a standard manner. We include a detailed proof here. For $k \in \bar{S}$, write

$$f(k, T) = \sum_{m \geq 0} f_m(k) T^m.$$

It suffices to prove that for each fixed $m$, the family $f_m(k)$ parametrized by $k \in \bar{S}$ is uniformly continuous. Our uniform continuity assumption means that for any $\varepsilon > 0$, there is a real number $\delta(\varepsilon) > 0$ such that for all $k_1$ and $k_2$ in $S$,

$$|f_m(k_1) - f_m(k_2)| < \varepsilon, \quad \text{whenever } d(k_1, k_2) < \delta(\varepsilon).$$

If we let $k_2$ vary, the continuity of $f_m(k_2)$ in $k_2$ shows that for all $k_1 \in S$ and $k_2 \in \bar{S}$,

$$|f_m(k_1) - f_m(k_2)| < \varepsilon, \quad \text{whenever } d(k_1, k_2) < \frac{1}{2}\delta(\varepsilon).$$

If we now let $k_1$ vary, the continuity of $f_m(k_1)$ in $k_1$ shows that for all $k_1 \in \bar{S}$ and $k_2 \in \bar{S}$,

$$|f_m(k_1) - f_m(k_2)| < \varepsilon, \quad \text{whenever } d(k_1, k_2) < \frac{1}{4}\delta(\varepsilon).$$

The corollary is proved. $\qquad\qquad\qquad\qquad\qquad\qquad\qquad\qquad\qquad\qquad\square$

To obtain similar result for meromorphic families, we need to introduce another notion about families of functions.

*Definition 5.3.* A family $f(k, T)$ is called a *strong family* in $|T|_\pi < p^c$ if $f(k, T)$ is the quotient of two uniformly continuous families of uniformly analytic functions in $|T|_\pi < p^c$.

There is probably a better choice of terminology for the notion of a strong family. We have tried to avoid the weaker notion of a uniformly continuous family of uniformly meromorphic functions, which simply means a uniformly continuous family which is also uniformly meromorphic. In such a family, it is not immediately obvious (although probably true) that the limit function is still a meromorphic function in the expected disc. In our stronger definition, this property is built in because of Corollary 5.2. Thus, we have:

COROLLARY 5.4.    *If $f(k, T)$ is a strong family in $|T|_\pi < p^c$ for $k \in S$, then $f(k, T)$ extends uniquely to a strong family in $|T|_\pi < p^c$ for $k \in \bar{S}$ (the closure of $S$).*

Sometimes, we shall need a slightly weaker notion of continuity. This will be needed later on. We give the definition here. A family of power series $f(k,T) \in R[[T]]$ is called *essentially continuous* in the disc $|T|_\pi < p^c$ if $f(k,T)$ can be written (in many different ways in general) as a product of two families

$$f(k,T) = f_1(k,T)f_2(k,T),$$

where $f_1(k,T)$ is a continuous family in $k$ and $f_2(k,T)$ (possibly not continuous in $k$) is a family of functions, analytic, without reciprocal zeros and with norm 1 in the disc $|T|_\pi \leq p^c$. Namely, $f_2(k,T)$ is a 1-unit in $|T|_\pi \leq p^c$. Two families of functions $f(k,T)$ and $g(k,T)$ are called *equivalent* in $|T|_\pi < p^c$ if $f(k,T) = g(k,T)h(k,T)$, where $h(k,T)$ is a 1-unit in $|T|_\pi \leq p^c$ for each $k \in S$. This is an equivalence relation. If one family $f(k,T)$ is essentially continuous and equivalent to another family $g(k,T)$ in $|T|_\pi < p^c$, then the second family $g(k,T)$ is also essentially continuous in $|T|_\pi < p^c$. A family $f(k,T)$ is said to be an essentially continuous family of uniformly analytic functions in the disc $|T|_\pi < p^c$ if it is both essentially continuous and uniformly analytic in the disc $|T|_\pi < p^c$. In terms of the above decomposition, this implies that the continuous part $f_1(k,T)$ is uniformly analytic in the disc $|T|_\pi < p^c$ because $1/f_2(k,T)$ is automatically uniformly analytic in that disc. If, in addition, the continuous part $f_1(k,T)$ is uniformly continuous, then we say that $f(k,T)$ is *essentially uniformly continuous*. A family $f(k,T)$ is called an *essentially strong family* in the disc $|T|_\pi < p^c$ if it is the quotient of two essentially uniformly continuous families of uniformly analytic functions in the disc $|T|_\pi < p^c$. In other words, an *essentially strong family* in the disc $|T|_\pi < p^c$ is a family which is equivalent to a strong family in the disc $|T|_\pi < p^c$. Thus, the notion of an essentially strong family in $|T|_\pi < p^c$ depends only on the equivalent class in $|T|_\pi < p^c$ of the family $f(k,T)$.

To say that a family $f(k,T)$ parametrized by $k \in S$ extends to an essentially strong family $g(k,T)$ parametrized by $k \in \bar{S}$ in $|T|_\pi < p^c$ means that $g(k,T)$ is an essentially strong family and its restriction to $k \in S$ is equivalent to $f(k,T)$. In particular, even for $k \in S$, $g(k,T)$ may be different from $f(k,T)$. But they are equivalent in $|T|_\pi < p^c$. Thus, in general, a family $f(k,T)$ parametrized by $k \in S$ may have many different ways to extend to an essentially strong family $g(k,T)$ parametrized by $k \in \bar{S}$ in $|T|_\pi < p^c$. By our definition, one sees that Corollary 5.4 extends to essentially continuous families.

COROLLARY 5.5.    *If $f(k,T)$ is an essentially strong family in $|T|_\pi < p^c$ for $k \in S$, then $f(k,T)$ extends (in many ways in general) to an essentially strong family $g(k,T)$ in $|T|_\pi < p^c$ for $k \in \bar{S}$ (the closure of $S$). Any two extended essentially strong families of $f(k,T)$ are equivalent in $|T|_\pi < p^c$.*

This result is true, because $f(k,T)$ is equivalent to a strong family in $|T|_\pi < p^c$ and any strong family extends uniquely. It is clear that finite products and quotients (if the denominator is not the zero family) of strong families (resp. essentially strong families) in $|T|_\pi < p^c$ for $k \in S$ are again strong families (resp. essentially strong families).

Our main concern will be the family of $L$-functions arising from a family of nuclear $\sigma$-modules. We want to know when such a family of $L$-functions is uniformly meromorphic and when it is a strong family. Before doing so, we need to define the notion of a uniform family of $\sigma$-modules. Let $B(k,X)$ be the matrix of a family of nuclear $\sigma$-modules $(M(k), \phi(k))$ parametrized by a parameter $k$. There may be no relations among the $\phi(k)$ for different $k$. In particular, the rank of $(M(k), \phi(k))$ could be totally different (some finite and others infinite, for instance) as $k$ varies. Write

$$B(k,X) = \sum_{u \in \mathbf{Z}^n} B_u(k) X^u,$$

where each coefficient matrix $B_u(k)$ has entries in $R$:

$$B_u(k) = (b_{w_1, w_2}(u, k)), \quad b_{w_1, w_2}(u, k) \in R,$$

where $w_1$ is the row index of $B_u(k)$ and $w_2$ is the column index of $B_u(k)$.

*Definition* 5.6. The family $B(k,X)$ (or the family $\phi(k)$) is called uniformly $c \log$-convergent if the following two conditions hold. First,

$$\lim_{|u| \to \infty} \inf \frac{\inf_k \operatorname{ord}_\pi B_u(k)}{\log_q |u|} \geq c.$$

Second, for any positive number $C > 0$, there is an integer $N_C > 0$ such that for all column numbers $w_2 > N_C$,

$$\operatorname{ord}_\pi b_{w_1, w_2}(u, k) \geq C,$$

uniformly for all $u, k, w_1$.

Note that the second condition above is automatically satisfied if the rank of $(M(k), \phi(k))$ is uniformly bounded. This is the case if each $\sigma$-module $M(k)$ has the same finite rank for every $k$. The following uniform result is of basic importance to our investigation.

THEOREM 5.7. *Let $B(k,X)$ be a family of matrices parametrized by $k$. Assume that the family $B(k,X)$ is uniformly $c \log$-convergent for some $0 < c < \infty$. Then for every $\varepsilon > 0$, the family of $L$-functions $L(B(k,X), T)$ is uniformly meromorphic in the open disc $|T|_\pi < p^{c-\varepsilon}$. If in addition, the family $L(B(k,X), T)$ is uniformly continuous (resp. essentially uniformly continuous) in the disc $|T|_\pi < p^{c-\varepsilon}$, then $L(B(k,X), T)$ is a strong family (resp. an essentially strong family) in the disc $|T|_\pi < p^{c-\varepsilon}$.*

A slightly weaker version of this uniform result is given in Theorem 5.2 of [26] in the case that the rank of the family $B(k, X)$ is a finite constant. We shall need the above more general result which allows the rank of $B(k, X)$ to be unbounded. We include a proof here, closely following the arguments in [26].

For each parameter $k$, let $F(k)$ be the following infinite matrix

$$F(k) = (B_{qu-v}(k))_{v,u \in \mathbf{Z}^n}$$

with block entries $B_{qu-v}(k)$, where $v$ denotes the row block index and $u$ denotes the column block index. Note that each block entry $B_{qu-v}(k)$ is a nuclear matrix. The Fredholm determinant $\det(I - TF(k))$ is defined and is $p$-adic meromorphic in $|T|_\pi < p^c$ since $B(k, X)$ is $c$ log-convergent. This is proved in [26] if $B(k, X)$ has finite rank. The infinite rank case is obtained in a similar way or simply by taking the limit. It also follows from the following stronger proof in Lemma 5.8. The Dwork trace formula for the $n$-torus is given by

$$L(B(k, X), T)^{(-1)^{n-1}} = \prod_{j=0}^{n} \det(I - q^j TF(k))^{(-1)^j \binom{n}{j}}.$$

This is known to be true if $B(k, X)$ has finite rank. The infinite rank case is similar. It can also be obtained from the finite rank case by taking the limit. The trace formula is a universal formula in some sense. Our two conditions on $B(k, X)$ guarantee that everything is convergent as shown below. Because of the shifting factor $q^j$, for the first part of the theorem, it suffices to prove the following lemma.

LEMMA 5.8.    *The Fredholm determinant* $\det(I - TF(k))$ *is a family of uniformly analytic functions in the finite open disc* $|T|_\pi < p^{c-\varepsilon}$.

*Proof.* For any $0 < \varepsilon < 2c$, we define a nonnegative weight function $w(u)$ on $\mathbf{Z}^n$ by

$$(5.2) \qquad w(u) = \begin{cases} (c - \frac{\varepsilon}{2}) \log_q |u|, & \text{if } |u| > 0, \\ 0, & \text{if } u = 0. \end{cases}$$

This function is increasing in $|u|$. Let $G_{v,u}(k)$ be the twisted matrix

$$G_{v,u}(k) = \pi^{w(v)-w(u)} B_{qu-v}(k).$$

This is also a nuclear matrix. The entries of this matrix will be in the quotient field of $R$, but they are bounded. Let $G(k)$ be the twisted infinite Frobenius matrix $(G_{v,u}(k))_{v,u \in \mathbf{Z}^n}$. Then,

$$\det(I - TF(k)) = \det(I - TG(k)).$$

Lemma 5.8 follows easily from the following lemma and a standard argument on determinant expansion; see the proof in [26].

LEMMA 5.9.    *For any $\varepsilon > 0$, all but finitely many column vectors $\vec{V}$ of $G(k)$ satisfy the inequality*

$$(5.3) \qquad\qquad \mathrm{ord}_\pi(\vec{V}) \geq c - \varepsilon.$$

*All (including those exceptional, finitely many, column vectors from (5.3)) of these are bounded by*

$$(5.4) \qquad\qquad \mathrm{ord}_\pi(\vec{V}) \geq -N(\varepsilon),$$

*where $N(\varepsilon)$ is a finite positive constant independent of $k$ and $\vec{V}$.*

*Proof.* For a given $\varepsilon > 0$ with $\varepsilon < 2c$, by our uniform $c$ log-convergent assumption, there is an integer $N_\varepsilon > 0$ such that for all $|u| > N_\varepsilon$,

$$(5.5) \qquad\qquad \mathrm{ord}_\pi B_u(k) \geq w(u) = (c - \frac{\varepsilon}{2}) \log_q |u|.$$

Take a positive integer $N_\varepsilon^*$ to be so large that $N_\varepsilon^* > N_\varepsilon$ and

$$(5.6) \qquad\qquad (c - \frac{\varepsilon}{2}) \log_q(q - \frac{N_\varepsilon}{N_\varepsilon^*}) \geq c - \varepsilon.$$

If $|u| \leq N_\varepsilon^*$, we have the trivial inequality $\mathrm{ord}_\pi G_{v,u}(k) \geq -w(u) \geq -w(N_\varepsilon^*)$. To prove the lemma, we first prove the claim that

$$\mathrm{ord}_\pi G_{v,u}(k) \geq c - \varepsilon$$

uniformly for all $|u| > N_\varepsilon^*$, all $v$ and all $k$. Assume $|u| > N_\varepsilon^*$. If $v = 0$ or $qu$, one checks that $\mathrm{ord}_\pi G_{v,u}(k) \geq w(qu) - w(u) = (c - \varepsilon/2) > c - \varepsilon$ and the claim is true. We now assume that $|u| > N_\varepsilon^*$ and $v$ is different from $0$ and $qu$. There are two cases.

If $|qu - v| \leq N_\varepsilon$, then by (5.2) and (5.6),

$$\begin{aligned}
\mathrm{ord}_\pi G_{v,u}(k) &\geq (c - \frac{\varepsilon}{2}) \log_q \frac{|v|}{|u|} \\
&= (c - \frac{\varepsilon}{2}) \log_q \frac{q|u| + (|v| - q|u|)}{|u|} \\
&\geq (c - \frac{\varepsilon}{2}) \log_q(q - \frac{N_\varepsilon}{N_\varepsilon^*}) \\
&\geq c - \varepsilon.
\end{aligned}$$

If $|qu - v| > N_\varepsilon$, we use the inequality $ab \geq a + b - 1$ (i.e., $(a-1)(b-1) \geq 0$ for $a, b \geq 1$) and deduce that

$$
\begin{aligned}
\mathrm{ord}_\pi G_{v,u}(k) &\geq (c - \frac{\varepsilon}{2}) \log_q \frac{|v||qu - v|}{|u|} \\
&\geq (c - \frac{\varepsilon}{2}) \log_q \frac{|v| + |qu - v| - 1}{|u|} \\
&\geq (c - \frac{\varepsilon}{2}) \log_q (\frac{q|u| - 1}{|u|}) \\
&\geq (c - \frac{\varepsilon}{2}) \log_q (q - \frac{1}{N_\varepsilon^*}) \\
&\geq c - \varepsilon.
\end{aligned}
$$

The claim is proved. It shows that all column vectors $\vec{V}$ contained in the block column of $G(k)$ indexed by $u$ with $|u| > N_\varepsilon^*$ satisfy the inequality in (5.3). To finish the proof of this inequality, we restrict our attention to the column vectors of $G(k)$ contained in the finitely many block columns of $G(k)$ indexed by $|u| \leq N_\varepsilon^*$. For each such $u$, our second condition in Definition 5.6 shows that only a uniformly bounded number of column vectors in $G(k)$ may not satisfy the inequality in (5.3). Since there are only finitely many such $u$, inequality (5.3) is proved. To get (5.4), we can assume that $|u| < N_\varepsilon^*$ by (5.3). For these small $u$, inequality (5.4) holds with $N(\varepsilon) = w(N_\varepsilon^*)$. The lemma is proved.

It remains to prove the second part of the theorem. By the inverted version of the Dwork trace formula, we have

$$
\det(I - TF(k)) = \prod_{j=0}^\infty L(B(k, X), q^j T)^{(-1)^{n-1}\binom{n+j-1}{j}}.
$$

Since the family $L(B(k, X), T)$ is uniformly continuous (resp. essentially uniformly continuous), the above shifting factor $q^j$ shows that the family $\det(I - TF(k))$ is also uniformly continuous (resp., essentially uniformly continuous). By Lemma 5.8, the family $\det(I - TF(k))$ is a uniformly continuous (resp., essentially uniformly continuous) family of uniformly analytic functions in $|T|_\pi < p^{c-\varepsilon}$. By Dwork's trace formula again, the family $L(B(k, X), T)$ is a strong (resp. essentially strong) family in $|T|_\pi < p^{c-\varepsilon}$. The theorem is proved.                                                                              □

In order to apply the above uniform result to the study of Dwork's unit root zeta functions, we need to understand the uniform convergent properties for the family $\mathrm{Sym}^k \phi$ of $\sigma$-modules arising from various symmetric powers of a given nuclear $\sigma$-module $(M, \phi)$ as $k$ varies. Such a family is apparently not uniformly $c$ log-convergent in general, even if the initial one $\phi$ is $c$ log-convergent. However, we have the following result, which is at the heart of our uniform proof later on.

Lemma 5.10. *Let $(M, \phi)$ be a $c$ log-convergent nuclear $\sigma$-module for some $0 < c < \infty$, ordinary at slope zero. Let $B$ be the matrix of $\phi$ under a basis which is ordinary at slope zero. Assume that $h_0 = 1$ and that the reduction of $B$ modulo $\pi$ is a constant matrix (i.e., each nonconstant term of $B$ involving the variable $X$ is divisible by $\pi$). Then, the family $\mathrm{Sym}^k B$ parametrized by positive integer $k$ is uniformly $(c - \varepsilon)$ log-convergent for every $\varepsilon > 0$.*

There are two conditions in our definition of uniform $c$ log-convergence. The proof of the above lemma naturally splits into two parts as well, corresponding to checking each one of the two conditions. We first prove that the second condition of uniformity is satisfied. Write

$$\mathrm{Sym}^k B = \sum_{u \in \mathbf{Z}^n} B_u(k) X^u,$$

where the coefficients $B_u(k)$ are matrices with entries in $R$. The rank of $B_u(k)$ is equal to the rank of $\mathrm{Sym}^k \phi$, which is unbounded in general as $k$ varies, even if the rank of $B$ is assumed to be finite. Write

$$B_u(k) = (b_{w_1, w_2}(u, k)),$$

where $w_1$ denotes the row index and $w_2$ denotes the column index of the matrix $B_u(k)$.

Lemma 5.11. *For any finite constant $C$, there is an integer $C_1 > 0$ such that for all $w_2 > C_1$, there exists the inequality*

$$\mathrm{ord}_\pi b_{w_1, w_2}(u, k) \geq C$$

*uniformly for all $u, k, w_1$.*

*Proof.* Let $\vec{e} = \{e_1, e_2, \cdots\}$ be a basis of $M$ ordinary at slope zero. We order the basis

$$\{e_{i_1} e_{i_2} \cdots e_{i_k} \mid 1 \leq i_1 \leq i_2 \leq \cdots \leq i_k\}$$

of the $k^{\mathrm{th}}$ symmetric product $\mathrm{Sym}^k M$ in some order compatible with the increasing size of $i_1 + \cdots + i_k$. For instance, we order

$$\{e_1^k, e_1^{k-1} e_2, e_1^{k-2} e_2^2, e_1^{k-1} e_3, e_1^{k-2} e_2 e_3, \cdots\}.$$

By definition of the map $\mathrm{Sym}^k \phi$,

$$(5.7) \qquad \mathrm{Sym}^k \phi(e_{i_1} \cdots e_{i_k}) = \phi(e_{i_1}) \phi(e_{i_2}) \cdots \phi(e_{i_k}).$$

Because $\phi$ is nuclear, the inequality

$$\mathrm{ord}_\pi \phi(e_i) \geq C$$

holds for all $i > r$, where $r$ is a finite number depending on the given integer $C$. Thus, in equation (5.7), we may restrict our attention to the indices in the range

$$1 \leq i_1 \leq i_2 \leq \cdots \leq i_k \leq r.$$

Since $h_0 = 1$ and our basis is ordinary at slope zero, $\phi(e_i)$ is divisible by $\pi$ for every $i \geq 2$. This shows that the product $\phi(e_{i_1})\phi(e_{i_2})\cdots\phi(e_{i_k})$ is divisible by $\pi^C$ if at least $C$ of the indices in $\{i_1, \cdots, i_k\}$ are greater than one. If $\{i_1, \cdots, i_k\}$ is a sequence with $1 \leq i_1 \leq i_2 \leq \cdots \leq i_k \leq r$ such that at most $C$ of the indices in it are greater than one, then we must have

$$1 = i_1 = i_2 = \cdots = i_{k-C} \leq i_{k-C+1} \leq \cdots \leq i_k \leq r.$$

The number of such sets $\{i_1, \cdots, i_k\}$ is at most $r^C$. Thus, for $w_2 \geq r^C + 1$, we have the inequality

$$\mathrm{ord}_\pi b_{w_1, w_2}(u, k) \geq C$$

uniformly for all $u, k, w_1$. The lemma is proved.                                       □

Next, we prove that the first condition of uniformity is satisfied.

LEMMA 5.12.    *For any $\varepsilon > 0$, there is a constant $N_\varepsilon > 0$ such that the inequality*

$$\mathrm{ord}_\pi B_u(k) \geq (c - \varepsilon)\log_q |u|$$

*holds uniformly for all $|u| > N_\varepsilon$ and all $k$.*

*Proof.* Let $(b_{ij}(X))$ be the matrix of $\phi$ under a basis $\{e_1, e_2, \cdots\}$ which is ordinary at slope zero. Then, we can write

$$\phi(e_j) = \sum_i b_{ij}(X)e_i = \sum_i (\sum_u b_{ij}^{(u)} X^u)e_i, \ b_{ij}(X) \in A_c.$$

Expanding the product on the right side of (5.7), we get

$$\phi(e_{j_1})\cdots\phi(e_{j_k}) = \sum_{i_1,\cdots,i_k} \sum_{u^{(1)},\cdots,u^{(k)} \in \mathbf{Z}^n} b_{i_1 j_1}^{(u^{(1)})} X^{u^{(1)}} \cdots b_{i_k j_k}^{(u^{(k)})} X^{u^{(k)}} e_{i_1} \cdots e_{i_k}.$$

We need to show that if $|u^{(1)}| + \cdots + |u^{(k)}| > N_\varepsilon$, then

$$\mathrm{ord}_\pi(b_{i_1 j_1}^{(u^{(1)})} \cdots b_{i_k j_k}^{(u^{(k)})}) \geq (c - \varepsilon)\log_q(|u^{(1)}| + \cdots + |u^{(k)}|).$$

By using a smaller positive integer $k$ if necessary, we may assume that all the exponents $u^{(1)}, \cdots, u^{(k)}$ in a typical term of the above expansion are nonzero. In this case, the coefficients $b_{i_\ell j_\ell}^{(u^{(\ell)})}$ are divisible by $\pi$ since $B$ is a constant matrix modulo $\pi$. Thus, we can write

$$b_{ij}^{(u)} = \pi a_{ij}^{(u)}$$

in $R$ for all $u \neq 0$, $i$ and $j$. The ring $A_c$ is $\pi$-saturated; i.e., if $f = \pi f_1$ for some $f_1 \in A_0$ and $f \in A_c$, then $f_1$ is also in $A_c$. We choose $N_\varepsilon$ sufficiently large such that the inequality

$$\mathrm{ord}_\pi a_{ij}^{(u)} \geq (c - \varepsilon) \log_q \frac{|u|}{N_\varepsilon}$$

holds for all $u \neq 0$, $i$ and $j$. This choice is possible for $b_{ij}^{(u)}$ and thus possible for $a_{ij}^{(u)}$ as well since

$$\mathrm{ord}_\pi a_{ij}^{(u)} = \mathrm{ord}_\pi b_{ij}^{(u)} - 1.$$

From these, we deduce that for all $u^{(\ell)} \neq 0$,

$$\begin{aligned}
\mathrm{ord}_\pi(b_{i_1 j_1}^{(u^{(1)})} \cdots b_{i_k j_k}^{(u^{(k)})}) &= k + \mathrm{ord}_\pi(a_{i_1 j_1}^{(u^{(1)})} \cdots a_{i_k j_k}^{(u^{(k)})}) \\
&\geq k + \max_{1 \leq \ell \leq k} \mathrm{ord}_\pi a_{i_\ell j_\ell}^{(u^{(\ell)})} \\
&\geq k + (c - \varepsilon) \log_q \frac{|u^{(1)}| + \cdots + |u^{(k)}|}{k N_\varepsilon} \\
&\geq k + (c - \varepsilon) \log_q(|u^{(1)}| + \cdots + |u^{(k)}|) - (c - \varepsilon) \log_q k N_\varepsilon \\
&\geq (c - 2\varepsilon) \log_q(|u^{(1)}| + \cdots + |u^{(k)}|)
\end{aligned}$$

uniformly for all sufficiently large $|u^{(1)}| + \cdots + |u^{(k)}|$ and all $k$. The lemma is proved. $\qquad\square$

To obtain more uniformly $c$ log-convergent families from known ones, one can use operations such as direct sum and tensor product. That is, the set of uniformly $c$ log-convergent families is stable under direct sum and tensor product. This is obvious for direct sum. We shall prove the stability for the tensor operation.

LEMMA 5.13. *Let $B(k, X)$ and $G(k, X)$ be two families of uniformly $c$ log-convergent matrices. Then their tensor product $B(k, X) \otimes G(k, X)$ is also a family of uniformly $c$ log-convergent matrices.*

*Proof.* Write

$$B(k, X) = \sum_{u \in \mathbf{Z}^n} B_u(k) X^u, \quad G(k, X) = \sum_{u \in \mathbf{Z}^n} G_u(k) X^u.$$

Then

$$B(k, X) \otimes G(k, X) = \sum_{u, v} B_u(k) \otimes G_v(k) X^{u+v}.$$

For any $\varepsilon > 0$, there is an integer $N_\varepsilon > 0$ such that for all $|u| > N_\varepsilon$ and all $k$,

$$\mathrm{ord}_\pi B_u(k) \geq (c - \varepsilon) \log_q |u|, \quad \mathrm{ord}_\pi G_u(k) \geq (c - \varepsilon) \log_q |u|$$

and
$$(c - \varepsilon) \log_q \frac{|u|}{2} \geq (c - 2\varepsilon) \log_q |u|.$$

Thus, for all $|u + v| > 2N_\varepsilon$ and for all $k$,
$$\begin{aligned}
\operatorname{ord}_\pi(B_u(k) \otimes G_v(k)) &\geq \max(\operatorname{ord}_\pi B_u(k), \operatorname{ord}_\pi G_v(k)) \\
&\geq (c - \varepsilon) \log_q \frac{|u| + |v|}{2} \\
&\geq (c - 2\varepsilon) \log_q |u + v|.
\end{aligned}$$

This verifies the first condition. It remains to check the second condition. Write
$$B_u(k) = (b_{w_1, w_2}(u, k)), \quad G_u(k) = (g_{w_1, w_2}(u, k))$$

where $w_1$ is the row index and $w_2$ is the column index. For any positive number $C > 0$, there is an integer $N_C > 0$ such that for all column indexes $w_2 > N_C$, we have
$$\operatorname{ord}_\pi b_{w_1, w_2}(u, k) \geq C, \quad \operatorname{ord}_\pi g_{w_1, w_2}(u, k) \geq C$$

uniformly for all $u, k, w_1$. This shows that except for the first $N_C^2$ columns, all column vectors $\vec{V}$ of $B_u(k) \otimes G_v(k)$ satisfy
$$\operatorname{ord}_\pi(\vec{V}) \geq C$$

uniformly for all $u$, $v$ and $k$. The lemma is proved.  $\square$

COROLLARY 5.14. *Let $\phi$ and $\psi$ be two $c$ log-convergent nuclear $\sigma$-modules, both ordinary at slope zero. Let $B$ (resp. $G$) be the matrix of $\phi$ (resp. $\psi$) under a basis ordinary at slope zero. Assume that $h_0(\phi) = h_0(\psi) = 1$. Assume further that the reductions of $B$ and $G$ modulo $\pi$ are constant matrices. Then, the family $\operatorname{Sym}^k B \otimes \operatorname{Sym}^k G$ parametrized by $k$ is uniformly $(c - \varepsilon)$ log-convergent for every $\varepsilon > 0$.*

COROLLARY 5.15. *Let $\phi(k)$ be a family of uniformly $c$ log-convergent nuclear $\sigma$-modules. Let $\psi$ be a fixed $c$ log-convergent $\sigma$-module. Then the twisted family $\phi(k) \otimes \psi$ is also uniformly $c$ log-convergent.*

These two corollaries can be combined to give many uniformly $c$ log-convergent families of nuclear $\sigma$-modules.

## 6. The unit root theorem

We are now ready to start the proof of the unit root theorem. In fact, we will prove something stronger. Namely, the unit root functions $L(\phi_0^k, T)$ are not just meromorphic, but in fact they form a strong family in the expected disc in the $p$-adic topology of $k$.

THEOREM 6.1.     *Let $(M, \phi)$ be a $c$ log-convergent nuclear $\sigma$-module $(M, \phi)$ for some $0 < c < \infty$, ordinary at slope zero. Let $(U_0, \phi_0)$ be the unit root part of $(M, \phi)$. Assume that the rank $h_0$ of $\phi_0$ is one. Then, the family $L(\phi_0^k, T)$ of unit root zeta functions parametrized by $k$ is a strong family in the open disc $|T|_\pi < p^{c-\varepsilon}$ for any $\varepsilon > 0$, where $k$ varies in any given residue class modulo $q - 1$ with induced $p$-adic topology.*

*Proof.* Under our assumption, the matrix of $\phi$ under an ordinary basis has the following form

$$(6.1) \qquad\qquad B = \begin{pmatrix} B_{00} & \pi B_{01} \\ B_{10} & \pi B_{11} \end{pmatrix},$$

where $B_{00}$ is an invertible element of $A_c$, and the other $B_{ij}$ are matrices over $A_c$. We first assume that the reduction of $B$ modulo $\pi$ is a constant matrix all of whose entries are zero except for the top left corner entry which is 1. Namely, we assume that

$$(6.2) \qquad\qquad B_{00} \equiv 1 \ (\mathrm{mod}\,\pi), \ \ \mathrm{ord}_\pi B_{10} \geq 1.$$

Under this simplifying condition, we show that the family $L(\phi_0^k, T)$ is a strong family, where $k$ varies in the ring of integers with induced $p$-adic topology.

Let $k$ be an integer and let

$$k_m = k + p^m.$$

Then as $p$-adic numbers, we have $\lim_{m \to \infty} k_m = k$. Furthermore, $k_m$ is positive for all large $m$. For any $p$-adic 1-unit $\alpha \in R$, one sees easily that

$$\lim_{m \to \infty} \alpha^{k_m} = \alpha^k.$$

Since $\phi_0$ is of rank 1 and $B_{00}$ is a 1-unit fibre by fibre, the characteristic root of each Euler factor in $L(\phi_0, T)$ is a $p$-adic 1-unit in $R$. Thus, we have the following limiting formula

$$\begin{aligned}
(6.3) \qquad L(\phi_0^k, T) &= \lim_{m \to \infty} L(\phi_0^{k_m}, T) \\
&= \lim_{m \to \infty} L(\phi^{k_m}, T).
\end{aligned}$$

The second equality holds because $\lim_m \pi^{k_m} = 0$. Let $S$ be the sequence $\{k_m\}$. There are two obvious topologies on $S$. One is the $p$-adic topology which has $k$ as the unique limiting point not in $S$. The other is the sequence topology which has $\infty$ as the unique limiting point not in $S$. It is easy to check that these two topologies agree on the set $S$, although the distance function may be a little different. Thus, one could use either topology in (6.3). We use the sequence topology to stress that $m$ will go to infinity. To show that the limit function $L(\phi_0^k, T)$ is meromorphic in the expected disc, we need to show that the family $L(\phi^{k_m}, T)$ is a strong family with respect to the sequence topology.

Since $\text{Sym}^{k-i}\phi = 0$ for $k < i$, the following fundamental decomposition formula

$$L(\phi^k, T) = \prod_{i \geq 1} L(\text{Sym}^{k-i}\phi \otimes \wedge^i \phi, T)^{(-1)^{i-1}i}$$

is a finite product. Thus, we need to show that for each fixed $i$, the family $L(\text{Sym}^{k_m-i}\phi \otimes \wedge^i \phi, T)$ is a strong family, where $k_m$ varies over the set of positive integers in $S$ with the sequence topology. We first prove the uniform continuity property of this family.

LEMMA 6.2.    *The family of L-functions $L(\text{Sym}^{k_m-i}\phi \otimes \wedge^i \phi, T)$ is uniformly continuous in $k_m$ with respect to the sequence topology.*

*Proof.* With respect to the sequence topology, the uniformly continuous notion is the same as the continuous notion. Thus, it suffices to prove that $L(\text{Sym}^{k_m-i}\phi \otimes \wedge^i \phi, T)$ is continuous in $k_m$ with respect to the sequence topology. Since there are only a finite number of closed points with a bounded degree, it suffices to check that each Euler factor is continuous in $k$. At a closed point $\bar{x} \in \mathbf{G}_m/\mathbf{F}_q$, let $\alpha_0, \alpha_1, \cdots$ denote the characteristic roots (counting multiplicities) of $\phi$ at the fibre $\bar{x}$ such that $\alpha_0$ is the $p$-adic 1-unit and all $\alpha_j$ are divisible by $\pi$ for $j \geq 1$. One checks that each $\alpha_j$ for $j \geq 1$ is in fact divisible by the higher power $\pi^{\deg(x)}$. The local Euler factor at $x$ of $L(\text{Sym}^k\phi \otimes \wedge^i \phi, T)$ is given by the reciprocal of the product

$$(6.4) \qquad E_k(x, T) = \prod(1 - \alpha_{j_1} \cdots \alpha_{j_k} \alpha_{\ell_1} \cdots \alpha_{\ell_i} T^{\deg(x)}),$$

where the product runs over all

$$0 \leq j_1 \leq \cdots \leq j_k, \quad 0 \leq \ell_1 < \cdots < \ell_i.$$

The above product is convergent as $\alpha_j$ is more and more divisible by $\pi$ as $j$ increases. Let $B_x$ denote the diagonal matrix $\text{diag}(\alpha_0, \alpha_1, \cdots)$ and let $B_{x,1}$ denote the diagonal sub-matrix $\text{diag}(\alpha_1, \alpha_2, \cdots)$. With $k$ replaced by $k_m - i$ in (6.4), one checks that

$$E_{k_m-i}(x, T) = \prod_{j=0}^{k_m-i} \det(I - T^{\deg(x)}\alpha_0^{k_m-i-j}\text{Sym}^j B_{x,1} \otimes \wedge^i B_x).$$

Now, since $\alpha_0$ is a 1-unit, for fixed $i$ and $j$, the sequence $\alpha_0^{k_m-i-j}$ converges to $\alpha_0^{k-i-j}$ as $m$ goes to infinity. Since $\text{Sym}^j B_{x,1}$ is divisible by $\pi^{j \deg(x)}$, we deduce that the sequence $L(\text{Sym}^{k_m-i}\phi \otimes \wedge^i \phi, T)$ is continuous in $k_m$ with respect to the sequence topology. The lemma is proved.    □

We now return to the proof of Theorem 6.1. Because of the special form of our matrix $B(X)$, Lemma 5.10 and Corollary 5.15 show that the family $\text{Sym}^{k_m-i}\phi \otimes \wedge^i \phi$ is uniformly $(c - \varepsilon)$ log-convergent for each fixed $i$. We just proved that the $L$-function of this sequence is uniformly continuous. Thus, by Theorem 5.7, the family of $L$-functions $L(\mathrm{Sym}^{k_m-i}\phi \otimes \wedge^i\phi, T)$ is a strong family in $|T|_\pi < p^{c-\varepsilon}$ for each fixed $i$. By Corollary 5.4 and the limiting formula

$$(6.5) \qquad L(\phi_0^k, T) = \lim_{m \to \infty} \prod_{i \geq 1} L(\mathrm{Sym}^{k_m-i}\phi \otimes \wedge^i\phi, T)^{(-1)^{i-1}i},$$

we conclude that the limiting unit root zeta function $L(\phi_0^k, T)$ is meromorphic in the expected disc for each fixed integer $k$. Here, we are using the fact that $\wedge^i\phi$ is more and more divisible by $\pi$ as $i$ grows. Also, we are using the sequence topology in the limiting formula. This finishes the proof of the meromorphy part (for each $k$) of Theorem 6.1.

Once we know that each $L(\phi_0^k, T)$ is meromorphic, we want to consider it as a family of functions parametrized by the integer $k \in \mathbf{Z}$ with induced $p$-adic topology. We want to show that this family is a strong family. There are at least two approaches to do this. Here we describe one of them. The other approach is given in Section 8.

The point is that the family $L(\mathrm{Sym}^k\phi \otimes \wedge^i\phi, T)$ (and also the family $L(\phi^k, T)$) is not continuous in $k$ with respect to the $p$-adic topology. However, it is essentially continuous if we restrict to a fixed finite disc and remove a few small values of $k$. The above proof then goes through if we use the essentially continuous family. As the final family $L(\phi_0^k, T)$ is indeed continuous with respect to the $p$-adic topology, we will be able to drop the word "essentially" in the final conclusion.

To carry out this idea, we need to modify Lemma 6.2 as follows.

LEMMA 6.3. *Let $S_c$ be the set of all positive integers $k \geq c$ with induced $p$-adic topology. Then the family of $L$-functions $L(\mathrm{Sym}^k\phi \otimes \wedge^i\phi, T)$ is essentially uniformly continuous in $k$ in the disc $|T|_\pi < p^c$.*

*Proof.* The proof is similar to the proof of Lemma 6.2. In the product decomposition

$$E_k(x, T) = \prod_{j=0}^k \det(I - T^{\deg(x)}\alpha_0^{k-j}\mathrm{Sym}^j B_{x,1} \otimes \wedge^i B_x),$$

we just take $j$ runs up to $c + 1$ instead of $k$. Namely, we let

$$E_{k,c}(x, T) = \prod_{0 \leq j < c+1} \det(I - T^{\deg(x)}\alpha_0^{k-j}\mathrm{Sym}^j B_{x,1} \otimes \wedge^i B_x).$$

The function $E_{k,c}(x, T)$ is clearly uniformly continuous in $k \in S_c$ with respect to the $p$-adic topology since each of its factors is uniformly continuous in $k$ and we have a bounded number of factors. The quotient $E_k(x, T)/E_{k,c}(x, T)$ has no contribution of zeros and poles in the disc $|T|_\pi \leq p^c$. This is because for $j \geq c + 1$, $\mathrm{Sym}^j B_{x,1}$ is divisible by $\pi^{j\deg(x)}$ and hence divisible by $\pi^{(c+1)\deg(x)}$. This shows that $E_k(x, T)$ and thus $L(\mathrm{Sym}^k\phi \otimes \wedge^i\phi, T)$ is essentially uniformly continuous in the disc $|T|_\pi < p^c$ with respect to the $p$-adic topology of $k \in S_c$. The lemma is proved.                                                                                          □

We now return to the proof of the strong family part of Theorem 6.1. Lemma 6.3 and Theorem 5.7 show that for each fixed $i$, the family

$$f(k, i, T) = L(\mathrm{Sym}^k\phi \otimes \wedge^i\phi, T)$$

parametrized by $k \in S_c$ is an essentially strong family in $|T|_\pi < p^{c-\varepsilon}$. For a positive integer $k \in S_c$, let

$$g(k, i, T) = \lim_{m \to \infty} L(\mathrm{Sym}^{k+p^m}\phi \otimes \wedge^i\phi, T).$$

Since

$$\lim_{m \to \infty}(k + p^m) = k,$$

the Euler factor proof in Lemma 6.3 shows that the quotient $f(k, i, T)/g(k, i, T)$ is a 1-unit in $|T|_\pi \leq p^c$. Note that $f(k, i, T)$ is not equal to $g(k, i, T)$ since the family $f(k, i, T)$ is not continuous but only essentially continuous in $|T|_\pi < p^c$. Thus, the family $g(k, i, T)$ parametrized by $k \in S_c$ is equivalent to the family $f(k, i, T)$ in $|T|_\pi < p^c$. This implies that $g(k, i, T)$ parametrized by $k \in S_c$ is also an essentially strong family in $|T|_\pi < p^{c-\varepsilon}$. By (6.5), we have the formula for $k \in S_{2c}$:

$$L(\phi_0^k, T) = h(k, T) \prod_{1 \leq i \leq c} g(k - i, i, T),$$

where $h(k, T)$ is the product of those factors with $i > c$ and $h(k, T)$ is a 1-unit in $|T|_\pi \leq p^c$. Thus, the unit root family $L(\phi_0^k, T)$ parametrized by $k \in S_{2c}$ is an essentially strong family in $|T|_\pi < p^{c-\varepsilon}$. In this way, we can write

$$L(\phi_0^k, T) = f_1(k, T)f_2(k, T),$$

where $f_1(k, T)$ is a strong family in $|T|_\pi < p^{c-\varepsilon}$ and $f_2(k, T)$ is a family of analytic functions but without reciprocal zeros and bounded by 1 on $|T|_\pi \leq p^{c-\varepsilon}$. Since both $L(\phi_0^k, T)$ and $f_1(k, T)$ are uniformly continuous in $k$, we deduce that the family $f_2(k, T)$ is a uniformly continuous family of uniformly analytic functions (without reciprocal zeros) in $|T|_\pi \leq p^{c-\varepsilon}$. It follows that $L(\phi_0^k, T)$ is a strong family parametrized by $k \in S_{2c}$. This family extends uniquely to a strong family parametrized by $k$ varying in the topological closure of $S_{2c}$ which is the whole space $\mathbf{Z}_p$ of $p$-adic integers. By uniqueness and the continuity of the family $L(\phi_0^k, T)$, this extended family agrees with $L(\phi_0^k, T)$ for all $k \in \mathbf{Z}_p$. This concludes our proof under the above simpler condition (6.2) on $B$.                                                                           □

We now explain how to reduce the general case to the above simpler situation. Consider the new basis (still ordinary at slope zero)

$$(e_1, e_2, \cdots) \begin{pmatrix} I_{00} & 0 \\ B_{10}B_{00}^{-1} & I_{11} \end{pmatrix} = \vec{e}E$$

of $M$, where $I_{00}$ is the rank 1 identity matrix and $I_{11}$ is an identity matrix. The transition matrix $E$ is $c$ log-convergent since $B_{10}B_{00}^{-1}$ is $c$ log-convergent. One calculates that the matrix of $\phi$ under $\vec{e}E$ is given by

$$E^{-1}BE^\sigma = \begin{pmatrix} B_{00} + \pi B_{01}(B_{10}B_{00}^{-1})^\sigma & \pi B_{01} \\ \pi(B_{11} - B_{10}B_{00}^{-1}B_{01})(B_{10}B_{00}^{-1})^\sigma & \pi(B_{11} - B_{10}B_{00}^{-1}B_{01}) \end{pmatrix}.$$

This matrix is still $c$ log-convergent since both $E$ and $E^\sigma$ are $c$ log-convergent. Let $f$ be the polynomial which is obtained from the top left entry of the above matrix by dropping all terms divisible by $\pi$. Since $\phi$ is ordinary at slope zero, the polynomial $f$ is a unit in the coordinate ring of the $n$-torus over $\mathbf{F}_q$. This means that $f$ must be a monomial whose coefficient is a $p$-adic unit. It is of course an invertible element in $A_c$. Pulling out the rank one factor $f$, we see that in the $c$ log-convergent category, the $\sigma$-module $(M, \phi)$ is the tensor product of a rank one unit root $\sigma$-module with matrix $f$ and an ordinary nuclear $\sigma$-module with a matrix satisfying the above simpler condition (6.2). Namely, we have

$$\phi = f \otimes \psi,$$

where $f$ is of rank 1 and $\psi$ has a matrix satisfying (6.2). At each fibre $\bar{x}$, we have

$$\lim_{m \to \infty} (f(x)f(x^q) \cdots f(x^{q^{\deg(\bar{x})-1}}))^{(q-1)p^m} = 1.$$

Using the fact that $f$ is of rank 1, by looking at the Euler factors, one easily checks that

$$\begin{aligned} L(\phi_0^k, T) &= \lim_{m \to \infty} L(f^k \otimes \psi^{k+(q-1)p^m}, T) \\ &= \lim_{m \to \infty} \prod_{i \geq 1} L(f^k \otimes \mathrm{Sym}^{k+(q-1)p^m-i}\psi \otimes \wedge^i \psi, T)^{(-1)^{i-1}i}. \end{aligned}$$

We are now in a similar situation as before, except that we now have an extra rank 1 twisting factor $f^k$. For a fixed $k$ (and fixed $i$, of course), the family $f^k \otimes \mathrm{Sym}^{k+(q-1)p^m-i}\psi \otimes \wedge^i \psi$ parametrized by $p^m$ is uniformly $(c - \varepsilon)$ log-convergent. For a fixed $k$, $L(f^k \otimes \mathrm{Sym}^{k+(q-1)p^m-i}\psi \otimes \wedge^i \psi, T)$ is also continuous with respect to the sequence topology of the sequence $p^m$. In the same way as before, we deduce that the limiting function $L(\phi_0^k, T)$ is meromorphic in the open disc $|T|_\pi < p^{c-\varepsilon}$ for each fixed $k$.

To finish the proof, we have to show that $L(\phi_0^k, T)$ $(k \in S)$ is a strong family in $k$ with respect to the $p$-adic topology, where $S$ is the integers in a fixed residue class modulo $(q-1)$. This is done by the following trick of change of basis. Let $r$ denote the smallest nonnegative residue of $S$ modulo $q-1$. We observe that if we change our basis $\vec{e}$ to $g\vec{e}$, where $g$ is an invertible element in $A_c$, then the matrix of $\phi$ under $g\vec{e}$ will be $g^{-1}Bg^\sigma = g^\sigma g^{-1}B$, which is still ordinary at slope zero and $c$ log-convergent. In the case that $g$ is a monomial with coefficient 1, the new basis is then $g^{q-1}B$ and thus we can remove the factor $g^{q-1}$ without changing the $L$-function. Because our $f$ is indeed a monomial, we can thus write $f = ag$, where $a$ is a $p$-adic unit in $R$ and $g$ is a monomial with coefficient 1. In this way, we can replace $f^k$ by $a^k g^r$ in the above limiting formula since we can drop any power of $g^{q-1}$. Thus, we have shown, with $k_m = k + (q-1)p^m$, that

$$L(\phi_0^k, T) = \prod_{i \geq 1} \lim_{m \to \infty} L(a^{k_m}g^r \otimes \mathrm{Sym}^{k_m-i}\psi \otimes \wedge^i\psi, T)^{(-1)^{i-1}i},$$

where we have replaced $a^k$ by $a^{k_m}$ from the previous formula. Now $r$ is fixed and $a$ is a constant. Corollary 5.15 shows that the family $a^k g^r \otimes \mathrm{Sym}^{k-i}\psi \otimes \wedge^i\psi$ parametrized by $k \in S_c$ is uniformly $(c-\varepsilon)$ log-convergent, where $S_c$ consists of the elements of $S$ which are greater than $c$. It is clear that $a^k$ is uniformly continuous in $k \in S$. One checks as before that the family of $L$-functions $L(a^k g^r \otimes \mathrm{Sym}^{k-i}\psi \otimes \wedge^i\psi, T)$ is essentially uniformly continuous with respect to the $p$-adic topology of $S_c$ in the disc $|T|_\pi < p^{c-\varepsilon}$ and thus it forms an essentially strong family in the same disc. This implies that $L(\phi_0^k, T)$ is an essentially strong family on $S_c$. Since $L(\phi_0^k, T)$ is also a uniformly continuous family, it must be a strong family on $S_c$. The topological closure of $S_c$ includes all $S$. We conclude that $L(\phi_0^k, T)$ is a strong family on all of $S$. The theorem is proved.

## 7. The higher slope theorem

We are now ready to prove the higher slope theorem. Again, we will prove a stronger family version. The idea here will be refined and extended in [29] to treat the higher rank case of Dwork's conjecture.

THEOREM 7.1.        Let $(M, \phi)$ be a $c$ log-convergent nuclear $\sigma$-module $(M, \phi)$ for some $0 < c < \infty$, ordinary up to slope $j$ for some integer $j \geq 0$. Let $(U_i, \phi_i)$ be the unit root $\sigma$-module coming from the slope $i$ part of $(M, \phi)$ for $0 \leq i \leq j$. Assume that $h_j = 1$. Then the family of unit root zeta functions $L(\phi_j^k, T)$ parametrized by $k$ is a strong family in the open disc $|T|_\pi < p^{c-\varepsilon}$ for every $\varepsilon > 0$, where $k$ varies in any given residue class modulo $q-1$ with induced $p$-adic topology.

*Proof.* The case for $j = 0$ has already been proved by the previous unit root theorem. We may thus assume that $j \geq 1$. It is clear that

$$\wedge^{h_0 + \cdots + ih_i} \phi \equiv 0 \pmod{\pi^{h_1 + \cdots + ih_i}}.$$

Let $\psi_i$ be the twist

(7.1) $$\psi_i = \pi^{-(h_1 + \cdots + ih_i)} \wedge^{h_0 + \cdots + ih_i} \phi.$$

By (3.1), this is a $c$ log-convergent $\sigma$-module with $h_0(\psi_i) = 1$ and ordinary at slope zero for all $0 \leq i \leq j$. $\qquad\square$

LEMMA 7.2. *For an integer $k > 0$, there exists the limiting formula*

$$L(\phi_j^k, T) = \lim_{m \to \infty} L(\psi_{j-1}^{-k+(q-1)p^{m+k}} \otimes \psi_j^{k+(q-1)p^m}, T),$$

*where $\psi_j^k$ is the $k^{\text{th}}$ power (not tensor power) defined as at the end of Section 2.*

*Proof.* Since there are only a finite number of closed points with a bounded degree, it suffices to prove the corresponding limiting formula for each Euler factor. At a closed point $\bar{x} \in \mathbf{G}_m/\mathbf{F}_q$, let $\pi^{ih_i \deg(x)} \alpha_i$ be the product of the $h_i$ characteristic roots of $\phi$ at the fibre $\bar{x}$ with slope $i$, where $0 \leq i \leq j$. Similarly, let $\gamma_i$ be the unique characteristic root of $\psi_i$ at the fibre $\bar{x}$ which is a $p$-adic unit, where $0 \leq i \leq j$. Then for each $i$ with $0 \leq i \leq j$, we have the relation

$$\alpha_0 \alpha_1 \cdots \alpha_i = \gamma_i.$$

In particular,

$$\gamma_{j-1} \alpha_j = \gamma_j.$$

One checks that the elements $\alpha_i$ and $\gamma_i$ for $0 \leq i \leq j$ are actually units in $R$. For each $i$, let $\beta_i$ run over the nonunit characteristic roots of $\psi_i$. One computes that the local Euler factor at $x$ of $L(\psi_{j-1}^{-k+(q-1)p^{m+k}} \otimes \psi_j^{k+(q-1)p^m}, T)$ is given by the reciprocal of the product

$$
\begin{aligned}
E_{m,k}(x, T) =\ & (1 - (\gamma_{j-1})^{(q-1)(p^{m+k}+p^m)} \alpha_j^{k+(q-1)p^m} T^{\deg(x)}) \\
& \times \prod_{\beta_j} (1 - (\gamma_{j-1})^{-k+(q-1)p^{m+k}} \beta_j^{k+(q-1)p^m} T^{\deg(x)}) \\
& \times \prod_{\beta_{j-1}} (1 - \beta_{j-1}^{-k+(q-1)p^{m+k}} \gamma_j^{k+(q-1)p^m} T^{\deg(x)}) \\
& \times \prod_{\beta_{j-1}, \beta_j} (1 - \beta_{j-1}^{-k+(q-1)p^{m+k}} \beta_j^{k+(q-1)p^m} T^{\deg(x)}).
\end{aligned}
$$

Now both $\beta_{j-1}$ and $\beta_j$ are divisible by $\pi$. Also both $\gamma_{j-1}^{q-1}$ and $\alpha_j^{q-1}$ are 1-unit in $R$. We deduce that

$$\lim_{m \to \infty} E_{m,k}(x, T) = (1 - \alpha_j^k T^{\deg(x)}).$$

This is the same as the reciprocal of the Euler factor of $L(\phi_j^k, T)$ at $x$. The desired limiting formula of the lemma follows.

As in the proof of the unit root theorem, we can write

$$(7.2) \qquad\qquad \psi_i = f_i \otimes \varphi_i,$$

where $f_i$ is an invertible monomial in $A_c$, $\varphi_i$ is a $c$ log-convergent $\sigma$-module with $h_0(\varphi_i) = 1$ and ordinary at slope zero. Furthermore, the matrix of $\varphi_i$ under an ordinary basis satisfies the condition in (6.2). Since the $f_i$ have rank one, we can pull them out of the above limiting formula and get

$$L(\phi_j^k, T) = \lim_{m \to \infty} L(\Phi_m(k), T),$$

where

$$\Phi_m(k) = f_{j-1}^{-k+(q-1)p^{m+k}} \otimes f_j^{k+(q-1)p^m} \otimes \varphi_{j-1}^{-k+(q-1)p^{m+k}} \otimes \varphi_j^{k+(q-1)p^m}.$$

At each fibre $\bar{x}$, we have

$$\lim_{m \to \infty} (f(x)f(x^q) \cdots f(x^{q^{\deg(\bar{x})-1}}))^{(q-1)p^m} = 1.$$

Using this fact and the decomposition formula derived from Lemma 4.5, we deduce that

$$L(\phi_j^k, T) = \lim_{m \to \infty} L((\frac{f_j}{f_{j-1}})^k \otimes \varphi_{j-1}^{-k+(q-1)p^{m+k}} \otimes \varphi_j^{k+(q-1)p^m}, T)$$

$$= \lim_{m \to \infty} \prod_{\ell_1, \ell_2 \geq 1} L((\frac{f_j}{f_{j-1}})^k \otimes \Psi_m(k, \ell_1, \ell_2), T)^{(-1)^{(\ell_1+\ell_2)}\ell_1\ell_2},$$

where $\Psi_m(k, \ell_1, \ell_2)$ is given by

$$\mathrm{Sym}^{-k+(q-1)p^{m+k}-\ell_1}\varphi_{j-1} \otimes \wedge^{\ell_1}\varphi_{j-1} \otimes \mathrm{Sym}^{k+(q-1)p^m-\ell_2}\varphi_j \otimes \wedge^{\ell_2}\varphi_j.$$

For fixed $k, \ell_1$ and $\ell_2$, Corollaries 5.14 and 5.15 show that the family $\Psi_m(k, \ell_1, \ell_2)$ parametrized by the sequence $\{m\}$ is uniformly $c$ log-convergent. As in Lemma 6.2, one checks that the family of $L$-functions $L((\frac{f_j}{f_{j-1}})^k \otimes \Psi_m(k, \ell_1, \ell_2), T)$ is (uniformly) continuous with respect to the sequence topology of $\{m\}$. It is therefore a strong family by Theorem 5.7. This shows that the limiting function $L(\phi_j^k, T)$ is meromorphic in the expected disc for each fixed integer $k > 0$.

Let $S$ be the set of integers in a fixed residue class modulo $q-1$. To show that $L(\phi_j^k, T)$ is a strong family parametrized by $S$ with respect to the $p$-adic topology, one can use the notion of an essentially strong family as in the case $i = 0$. Let $S_c$ be the set consisting of the integers in $S$ which are greater than $c$. By the proof of Lemma 7.2, we also have the slightly different limiting formula:

$$L(\phi_j^k, T) = \lim_{a \to \infty} L(\psi_{j-1}^{-k_a+(q-1)p^{k_a}} \otimes \psi_j^{k_a}, T),$$

where $k_a$ is any sequence of positive integers in $S_c$ such that $\lim_a k_a = \infty$ as integers and $\lim_a k_a = k$ as $p$-adic integers. For instance, we can take $k_a = k + (q-1)p^a$. Using the decomposition formula from Lemma 4.5 and (7.2), we get

$$L(\phi_j^k, T) = \lim_a L(((\frac{f_j}{f_{j-1}})^{k_a} \otimes \varphi_{j-1}^{-k_a+(q-1)p^{k_a}} \otimes \varphi_j^{k_a}, T)$$

$$= \lim_a \prod_{\ell_1, \ell_2 \geq 1} L(((\frac{f_j}{f_{j-1}})^{k_a} \otimes \Psi(k_a, \ell_1, \ell_2), T)^{(-1)^{(\ell_1+\ell_2)}\ell_1\ell_2},$$

where

$$\Psi(k_a, \ell_1, \ell_2) = \operatorname{Sym}^{-k_a+(q-1)p^{k_a}-\ell_1}\varphi_{j-1} \otimes \wedge^{\ell_1}\varphi_{j-1} \otimes \operatorname{Sym}^{k_a-\ell_2}\varphi_j \otimes \wedge^{\ell_2}\varphi_j.$$

As in Lemma 6.3, one finds that the family $L(((\frac{f_j}{f_{j-1}})^{k_a} \otimes \Psi(k_a, \ell_1, \ell_2), T)$ of functions parametrized by $k_a \in S_c$ is essentially uniformly continuous in the expected disc. The family $\Psi(k_a, \ell_1, \ell_2)$ parametrized by $k_a$ is uniformly $c\log$-convergent by Corollaries 5.14 and 5.15. Using the trick of the change of basis, we can reduce the power of $(f_j/f_{j-1})^{k_a}$ modulo $q-1$ up to a twisting constant in $R$. This implies that the family $L(((\frac{f_j}{f_{j-1}})^{k_a} \otimes \Psi(k_a, \ell_1, \ell_2), T)$ parametrized by $k_a \in S_c$ is an essentially strong family. It extends to an essentially strong family on the topological closure of $S_c$. We deduce that the limiting family $L(\phi_j^k, T)$ is an essentially strong family on $S_c$ and hence on all of $S$. Because $L(\phi_j^k, T)$ is uniformly continuous in $k \in S$, we conclude that $L(\phi_i^k, T)$ is a strong family for $k \in S$ in the expected disc. The proof is complete. $\qquad\square$

## 8. The limiting $\sigma$-module and explicit formula

In the proof of the unit root theorem and the higher slope theorem, we used the notion of essential continuity to understand the variation of the unit root zeta function as $k$ varies with respect to the $p$-adic topology. There is another approach which avoids the essential continuity notion and works directly with continuous families. This will be useful for further investigations as well as some other applications. The purpose of this section is to describe this second approach in the optimal $c\log$-convergent setup. A full treatment of this approach in the simpler overconvergent setting will be given in a future paper when we need to get explicit information about zeros and poles of unit root zeta functions.

We assume that $\phi$ satisfies the simpler condition in (6.2) and that $k_m = k + p^m$. As we have seen before, the key is to understand the limit function of the sequence of $L$-functions $L(\operatorname{Sym}^{k_m}\phi, T)$ as $k_m$ varies. In the new approach presented in this section, we show that there is a $c\log$-convergent nuclear $\sigma$-module $(M_{\infty,k}, \phi_{\infty,k})$ called the limiting $\sigma$-module of the sequence $(\mathrm{Sym}^{k_m} M, \mathrm{Sym}^{k_m} \phi)$ such that

$$\lim_{m \to \infty} L(\mathrm{Sym}^{k_m} \phi, T) = L(\phi_{\infty,k}, T). \tag{8.1}$$

Furthermore, the family $(M_{\infty,k}, \phi_{\infty,k})$ parametrized by $k$ is uniformly $c\log$-convergent. As the positive integer $k$ varies $p$-adically, the family of functions $L(\phi_{\infty,k}, T)$ parametrized by $k$ turns out to be uniformly continuous in $k$. This provides another proof of the unit root theorem and the higher slope theorem. It shows that Dwork's unit root zeta function in the ordinary rank one case is an infinite alternating product of $L$-functions of nuclear $\sigma$-modules with the same convergent condition as the initial $\phi$. This infinite product can be taken to be a finite product if one restricts to a finite disc. In the case that the initial nuclear $\sigma$-module $(M, \phi)$ is of finite rank, the infinite product is actually a finite product. This can be viewed as a structure theorem. It implies that if the initial $(M, \phi)$ is $c\log$-convergent (resp. overconvergent), then Dwork's unit root zeta function in the ordinary rank one case can be expressed in terms of $L$-functions of $c\log$-convergent (resp. overconvergent) nuclear $\sigma$-modules of infinite rank.

Let $(M, \phi)$ be a nuclear $c\log$-convergent $\sigma$-module satisfying the simple condition in (6.2). Let $k$ be a positive integer. We want to identify the $k^{\text{th}}$ symmetric power $(\mathrm{Sym}^k M, \mathrm{Sym}^k \phi)$ with another $\sigma$-module $(M_k, \phi_k)$ where it is easier to take the limit as $k$ varies $p$-adically. Let $\vec{e} = \{e_1, e_2, \cdots\}$ be a basis of $(M, \phi)$ ordinary at slope zero such that its matrix satisfies the condition in (6.2). Define $M_k$ to be the Banach module over $A_0$ with the basis $\vec{f}(k)$:

$$\{f_{i_1} f_{i_2} \cdots f_{i_r} | 2 \le i_1 \le i_2 \le \cdots \le i_r, \ 0 \le r \le k\},$$

where we think of 1 (corresponding to the case $r = 0$) as the first basis element of $\vec{f}(k)$. There is an isomorphism of Banach $A_0$-modules between $\mathrm{Sym}^k M$ and $M_k$. This isomorphism is given by the map:

$$\Upsilon : e_1^{k-r} e_{i_1} \cdots e_{i_r} \longrightarrow f_{i_1} \cdots f_{i_r}.$$

Thus,

$$\Upsilon(e_1) = 1, \ \Upsilon(e_i) = f_i \ (i > 1).$$

We can give a nuclear $\sigma$-module structure on $M_k$. The semi-linear map $\phi_k$ acting on $M_k$ is given by the pull back $\Upsilon \circ \mathrm{Sym}^k \phi \circ \Upsilon^{-1}$ of $\mathrm{Sym}^k \phi$ acting on $\mathrm{Sym}^k M$. Namely,

$$\phi_k(f_{i_1} \cdots f_{i_r}) = \Upsilon(\phi(e_1^{k-r}) \phi(e_{i_1}) \cdots \phi(e_{i_r})). \tag{8.2}$$

Under the identification of $\Upsilon$, the map $\mathrm{Sym}^k \phi$ becomes $\phi_k$. Thus, the $k^{\text{th}}$ symmetric power $(\mathrm{Sym}^k M, \mathrm{Sym}^k \phi)$ is identified with $(M_k, \phi_k)$. It follows that

$$L(\mathrm{Sym}^k \phi, T) = L(\phi_k, T).$$

In order to take the limit of the sequence of modules $M_k$, we define $M_\infty$ to be the Banach module over $A_0$ with the countable basis $\vec{f}$:

$$(8.3) \qquad \{f_{i_1} f_{i_2} \cdots f_{i_r} | 2 \leq i_1 \leq i_2 \leq \cdots \leq i_r, \ 0 \leq r\},$$

where, again, we think of 1 (corresponding to the case $r = 0$) as the first basis element of $\vec{f}$. Note that in (8.3), there is no longer restriction on the size of $r$.

The module $M_\infty$ can be identified with the graded commutative algebra over $A_0$ generated by the free elements $\{f_2, f_3, \cdots\}$ which are of degree 1. Namely, as the graded commutative $A_0$-algebra,

$$M_\infty = A_0[[f_2, f_3, \cdots]],$$

where $A_0[[f]]$ denotes the ring of power series over $A_0$ in the infinite set of variables $\{f_2, f_3, \cdots\}$. To be more precise, we define the weight $w(f_i) = i$. Denote by $I$ the countable set of integer vectors $u = (u_2, u_3, \cdots)$ satisfying $u_i \geq 0$ and almost all of the $u_i$ (except for finitely many of them) are zero. For a vector $u \in I$, define

$$f^u = f_2^{u_2} f_3^{u_3} \cdots f_m^{u_m} \cdots, \quad w(u) = 2u_2 + 3u_3 + \cdots + mu_m + \cdots.$$

Then the infinite dimensional $A_0$-algebra $M_\infty$ can be described precisely as follows:

$$M_\infty = \{\sum_{u \in I} a_u f^u | a_u \in A_0\}.$$

The algebra $M_\infty$ becomes a Banach algebra over $A_0$ with a basis given by (8.3) or equivalently by $\{f^u | u \in I\}$. We order the elements of $I$ in any way which is compatible with the increasing size of $w(u)$. If $(M, \phi)$ is of finite rank, say of rank $r$, then one checks that the algebra $M_\infty$ is just the formal power series ring $A_0[[f_2, \cdots, f_r]]$ over $A_0$. Thus, the algebra $M_\infty$ is Noetherian if and only if the rank $r$ of $\phi$ is finite. Note that the Banach $A_0$-module $M_\infty$ has a countable basis but does not have an orthonormal basis.

For $u \in I$, the degree of $f^u$ is defined to be $\|u\| = u_2 + u_3 + \cdots$. The degree of an element $\sum a_u f^u$ in $M_\infty$ is defined to be the largest degree of its nonzero terms. Namely, the degree is given by

$$\sup\{\|u\| \mid a_u \neq 0\}.$$

The degree of an element $\sum a_u f^u$ could be infinite in general because there are infinitely many terms in $\sum a_u f^u$. For each positive integer $k$, the Banach module $M_k$ is the submodule of $M_\infty$ consisting of those elements of degree at most $k$. That is,

$$M_k = \{f \in M_\infty | \deg(f) \leq k\}.$$

For each integer $k$, define

$$M_{\infty,k} = M_\infty.$$

This module is thus independent of $k$. To define the nuclear semi-linear map $\phi_{\infty,k}$ on $M_{\infty,k}$ for each integer $k$, we take a sequence of positive integers $k_m$ (for instance $k_m = k + p^m$) such that $\lim k_m = \infty$ as integers and $\lim_m k_m = k$ as $p$-adic integers. Define $\phi_{\infty,k}$ by the following limiting formula:

$$(8.4) \qquad \phi_{\infty,k}(f_{i_1} \cdots f_{i_r}) = \lim_{m \to \infty} \Upsilon(\phi(e_1^{k_m - r})\phi(e_{i_1}) \cdots \phi(e_{i_r}))$$

$$= \big( \lim_{m \to \infty} \Upsilon(\phi(e_1^{k_m - r}))\big) \Upsilon(\phi(e_{i_1}) \cdots \phi(e_{i_r})).$$

To show that this is well-defined, we need to show that the above first limiting factor exists. Write

$$(8.5) \qquad\qquad\qquad \phi(e_1) = u(e_1 + \pi e),$$

where $u$ is a 1-unit in $R$ and $e$ is an element of $M$ which is an infinite linear combination of the elements in $\{e_2, e_3, \cdots\}$. By the binomial theorem, we have

$$\phi(e_1^k) = u^k(e_1 + \pi e)^k$$

$$= u^k(e_1^k + \binom{k}{1}\pi e_1^{k-1}e + \binom{k}{2}\pi^2 e_1^{k-2}e^2 + \cdots).$$

This identity implies that

$$(8.6)$$
$$\lim_m \Upsilon(\phi(e_1^{k_m - r})) = u^{k-r}(1 + \binom{k-r}{1}\pi\Upsilon(e) + \binom{k-r}{2}\pi^2\Upsilon(e)^2 + \cdots)$$

$$= u^{k-r}(1 + \pi\Upsilon(e))^{k-r}.$$

Thus, the map $\phi_{\infty,k}$ is well-defined for every integer $k$ and it is independent of the choice of our chosen sequence $k_m$. Furthermore, using (8.4) and (8.6), we see that $\phi_{\infty,k}$ makes sense for any $p$-adic integer $k$, not necessarily usual integers. With this definition, we have the following result.

THEOREM 8.1.    Let $(M, \phi)$ be a nuclear $c$ log-convergent $\sigma$-module which is ordinary at slope zero and where $h_0(\phi) = 1$. Assume that $\phi$ satisfies the simpler condition in (6.2). Then, the family $(M_\infty, \phi_{\infty,k})$ parametrized by the $p$-adic integer $k$ is a family of uniformly $c$ log-convergent nuclear $\sigma$-modules. Furthermore, the family of $L$-functions $L(\phi_{\infty,k}, T)$ parametrized by $p$-adic integer $k \in \mathbf{Z}_p$ is a strong family in the open disc $|T|_\pi < p^{c-\varepsilon}$ for any $\varepsilon > 0$.

*Proof.* The proof of the uniform part is essentially the same as the proof of Lemma 5.10. The main point is that all bounds in that proof are uniformly independent of $k$. We omit the details. For the strong family part, by Theorem 5.7, we only need to check that the family of $L$-functions $L(\phi_{\infty,k}, T)$ is uniformly continuous in $k$ with respect to the $p$-adic topology. By our definition of $\phi_{\infty,k}$ in (8.4) and (8.6), it suffices to show that

$$(8.7) \qquad \|u^{k_1-r}(1+\pi\Upsilon(e))^{k_1-r} - u^{k_2-r}(1+\pi\Upsilon(e))^{k_2-r}\| \le c\|k_1 - k_2\|$$

uniformly for all integers $r$, where $c$ is some positive constant. But this follows from the binomial theorem since $u$ is a 1-unit in $R$. The proof is complete. $\quad\square$

In order to take the limit of $\phi_k$ acting on the submodule $M_k$ of $M_\infty$ as $k$ varies in a sequence of positive integers, we semi-linearly extend $\phi_k$ to act on $M_\infty$ by requiring

$$\phi(f_{i_1} \cdots f_{i_r}) = 0$$

for all indices $2 \le i_1 \le i_2 \le \cdots \le i_r$ with $i_r > k$. In this way, $\phi_k$ induces a semi-linear endomorphism of the Banach algebra $M_\infty$ over $A_0$. The pair $(M_\infty, \phi_k)$ then becomes a nuclear $\sigma$-module over $A_0$. The $L$-function $L(\phi_k, T)$ is the same, whether we view $\phi_k$ as acting on $M_\infty$ or on the submodule $M_k$. Let $\mu(k)$ be the greatest integer such that for all integers $j > k$, we have

$$\phi(e_j) \equiv 0 \pmod{\pi^{\mu(k)}}.$$

The integers $\mu(k)$ are related to but different from the integers $d_k$ defined in Definition 2.7. Since $\phi$ is nuclear, it follows that

$$\lim_{k\to\infty} \mu(k) = \infty.$$

From our definitions of $\phi_k$ and $\phi_{\infty,k}$, one deduces that the following congruence formula holds.

LEMMA 8.2. *Let $k$ be an integer and let $k_m = k + p^m > 0$. Then as an endomorphism on $M_\infty$, there exists the congruence*

$$\phi_{k_m} \equiv \phi_{\infty,k} \pmod{\pi^{\min(\mu(k_m),m+1)}}.$$

*Proof.* We need to show that for all indices $2 \le i_1 \le i_2 \le \cdots \le i_r$,

$$\phi_{k_m}(f_{i_1} \cdots f_{i_r}) \equiv \phi_{\infty,k}(f_{i_1} \cdots f_{i_r}) \pmod{\pi^{\min(\mu(k_m),m+1)}}.$$

If $i_r > k_m$, the left side is zero and the right side is divisible by $\pi^{\mu(k_m)}$. The congruence is indeed true. Assume now that $i_r \le k_m$. By (8.2) and (8.4), it suffices to check that

$$(8.8) \qquad \Upsilon(\phi(e_1^{k_m-r})) \equiv \lim_a \Upsilon(\phi(e_1^{k_a-r})) \pmod{\pi^{m+1}},$$

where $k_a$ is a sequence of integers with $k$ as its $p$-adic limit and with $\infty$ as its complex limit. But (8.8) is a consequence of (8.5) and the binomial theorem. In fact, the congruence in (8.8) can be improved to

$$(8.9) \qquad \Upsilon(\phi(e_1^{k_m-r})) \equiv \lim_a \Upsilon(\phi(e_1^{k_a-r})) \pmod{\pi p^m}$$

if the ramification index $\mathrm{ord}_\pi p$ is smaller than $p-1$. The lemma is proved. $\square$

COROLLARY 8.3.    *Let $k$ be an integer and let $k_m = k + p^m$. Then, there is the limiting formula*

$$\lim_{m \to \infty} L(\mathrm{Sym}^{k_m}\phi, T) = \lim_{m \to \infty} L(\phi_{k_m}, T)$$
$$= L(\phi_{\infty,k}, T).$$

Using these results, we obtain an explicit formula for Dwork's unit root zeta function.

THEOREM 8.4.    *Let $(M, \phi)$ be a $c$ log-convergent nuclear $\sigma$-module for some $0 < c < \infty$, ordinary at slope zero with $h_0 = 1$. Let $(U_0, \phi_0)$ be the unit root part of $(M, \phi)$. Write $\phi = ag \otimes \psi$, where $a$ is a $p$-adic unit in $R$, $g$ is a monomial (hence in $A_c$) with coefficient 1 and $\psi$ satisfies the simpler condition in (6.2). Denote by $\psi_{\infty,k}$ the limiting nuclear $\sigma$-module of the sequence $\mathrm{Sym}^{k_m}\psi$. Then for all integers $k$ in the residue class of $r$ modulo $(q-1)$, there is the following explicit formula for Dwork's unit root zeta function:*

$$L(\phi_0^k, T) = \prod_{i \geq 1} L(a^k g^r \otimes \psi_{\infty,k-i} \otimes \wedge^i \psi, T)^{(-1)^{i-1}i}.$$

*In particular, the family $L(\phi_0^k, T)$ of $L$-functions parametrized by $k$ in any residue class modulo $(q-1)$ is a strong family in the open disc $|T|_\pi < p^{c-\varepsilon}$ with respect to the $p$-adic topology of $k$.*

It should be noted that this product is a finite product if $\phi$ is of finite rank, since we have $\wedge^i\psi = 0$ for $i$ greater than the rank of $\phi$. A similar formula holds for higher slope unit root zeta functions. We now make this precise. Its proof follows from the proof of Theorem 7.1.

THEOREM 8.5.    *Let $(M, \phi)$ be a $c$ log-convergent nuclear $\sigma$-module for some $0 < c < \infty$, ordinary up to slope $j$ for some integer $j \geq 0$. Let $(U_i, \phi_i)$ be the unit root $\sigma$-module coming from the slope $i$ part of $(M, \phi)$ for $0 \leq i \leq j$. Assume $h_j = 1$. Let $\psi_i$ be defined as in equation (7.1) for $0 \leq i \leq j$. Write $\psi_i = a_i g_i \otimes \varphi(i)$, where $a_i$ is a $p$-adic unit in $R$, $g_i$ is a monomial in $A_c$ with coefficient 1 and $\varphi(i)$ satisfies the simpler condition in (6.2). Denote by $\varphi_{\infty,k}(i)$ the limiting nuclear $\sigma$-module of the sequence $\mathrm{Sym}^{k_m}\varphi(i)$ as $m$ goes to infinity. Then for all integers $k$ in the residue class of $r$ modulo $(q-1)$, there is the following explicit formula for Dwork's $j$th unit root zeta function:*

$$L(\phi_j^k, T) = \prod_{\ell_1, \ell_2 \geq 1} L((\frac{a_j}{a_{j-1}})^k (\frac{g_j}{g_{j-1}})^r \otimes \Psi(k, \ell_1, \ell_2), T)^{(-1)^{(\ell_1+\ell_2)}\ell_1\ell_2},$$

*where*

$$\Psi(k, \ell_1, \ell_2) = \varphi_{\infty, -k - \ell_1}(j - 1) \otimes \wedge^{\ell_1} \varphi_{j-1} \otimes \varphi_{\infty, k - \ell_2}(j) \otimes \wedge^{\ell_2} \varphi_j.$$

*In particular, the family $L(\phi_j^k, T)$ of L-functions parametrized by $k$ in any residue class modulo $(q - 1)$ is a strong family in the open disc $|T|_\pi < p^{c - \varepsilon}$ with respect to the p-adic topology of $k$.*

In both Theorem 8.4 and Theorem 8.5, we could have used the decomposition formula derived from Lemma 4.8. If we did so, we would get explicit formulas with fewer cancellations but they would look more complicated (less compact). Thus, we shall omit them and use the above simpler looking formulas. For future applications to the higher rank case of Dwork's conjecture, we need a slightly more general result than Theorem 8.4. Namely, we need to twist the unit root family $\phi_0^k$ by a fixed nuclear $\sigma$-module $\varphi$. We state this generalization here. The proof is the same as the proof of Theorem 8.4. One simply twists the whole proof by $\varphi$. The twisted basic decomposition formula is

$$L(\phi^k \otimes \varphi, T) = \prod_{i=1}^{\infty} L(\mathrm{Sym}^{k-i} \phi \otimes \wedge^i \phi \otimes \varphi, T)^{(-1)^{i-1} i}.$$

The required uniform result is stated in Corollary 5.15.

THEOREM 8.6. *Let $(M, \phi)$ and $(N, \varphi)$ be two $c$ log-convergent nuclear $\sigma$-modules for some $0 < c < \infty$. Assume that the first one $\phi$ is ordinary at slope zero with $h_0 = 1$. Let $(U_0, \phi_0)$ be the unit root part of $(M, \phi)$. Write $\phi = ag \otimes \psi$, where $a$ is a p-adic unit in $R$, $g$ is a monomial in $A_c$ with coefficient $1$ and $\psi$ satisfies the simpler condition in (6.2). Denote by $\psi_{\infty, k}$ the limiting nuclear $\sigma$-module of the sequence $\mathrm{Sym}^{k_m} \psi$. Then for all integers $k$ in the residue class of $r$ modulo $(q - 1)$, there is the following explicit formula for Dwork's twisted unit root zeta function:*

$$L(\phi_0^k \otimes \varphi, T) = \prod_{i \geq 1} L(a^k g^r \otimes \psi_{\infty, k-i} \otimes \wedge^i \psi \otimes \varphi, T)^{(-1)^{i-1} i}.$$

*In particular, the twisted family $L(\phi_0^k \otimes \varphi, T)$ of L-functions parametrized by $k$ in the residue class of $r$ modulo $(q - 1)$ is a strong family of meromorphic functions in the disc $|T|_\pi < p^c$ with respect to the p-adic topology of $k$.*

If $\phi_0$ has rank greater than one, our current limiting method does not work. However, the embedding method introduced in [29] shows that there is always an infinite product formula for $L(\phi_0^k \otimes \varphi, T)$ even if $\phi_0$ has rank greater than 1. Such an infinite product formula comes from a totally different reason and cannot be written as a finite product when $\phi_0$ has rank greater than one even if $\phi$ itself is of finite rank; see [29].

University of California, Irvine, CA
*E-mail address*: dwan@math.uci.edu

(Received November 19, 1997)

(Revised November 20, 1998)